# Euler-Hurwitz series and non-linear Euler sums

Donal F. Connon

11 March 2008

**Abstract**

In this paper we derive two expressions for the Hurwitz zeta function $\varsigma(q+1,x)$ involving the complete Bell polynomials $Y_n(x_1,...,x_n)$ in the restricted case where $q$ is a positive integer and $q \geq 1$, one of which is

$$\varsigma(q+1,x) = \frac{1}{q!}\sum_{n=1}^{\infty}\frac{1}{n}\frac{\Gamma(n)\Gamma(x)}{\Gamma(n+x)}Y_{q-1}\left(0!H_n^{(1)}(x),\, 1!H_n^{(2)}(x),\, ...,\, (q-2)!H_n^{(q-1)}(x)\right)$$

where $H_n^{(m)}(x)$ are the generalised harmonic number functions defined by

$$H_n^{(m)}(x) = \sum_{k=0}^{n-1}\frac{1}{(k+x)^m}.$$

This in turn gives rise to "Euler-Hurwitz" series of the form

$$\varsigma(5,x) = \frac{1}{4!}\sum_{n=1}^{\infty}\frac{1}{n}\frac{\Gamma(n)\Gamma(x)}{\Gamma(n+x)}\left(\left[H_n^{(1)}(x)\right]^3 + 3H_n^{(2)}(x)H_n^{(1)}(x) + 2H_n^{(3)}(x)\right)$$

and similar expressions may then be used to determine identities for combinations of both linear and non-linear Euler sums; two such examples of which are shown below

$$4!\varsigma(5) = \sum_{n=1}^{\infty}\frac{\left[H_n^{(1)}\right]^3 + 3H_n^{(2)}H_n^{(1)} + 2H_n^{(3)}}{n^2}$$

$$\frac{1}{2}5!\varsigma(6) = \sum_{n=1}^{\infty}\frac{\left[H_n^{(1)}\right]^4 + 3\left[H_n^{(1)}\right]^2 H_n^{(2)} + H_n^{(1)}H_n^{(3)}}{n^2} - \sum_{n=1}^{\infty}\frac{\left[H_n^{(1)}\right]^3 + 3H_n^{(1)}H_n^{(2)} + H_n^{(3)}}{n^3}$$

INTRODUCTION

The Riemann zeta function $\varsigma(s)$ is defined for complex values $s$ by [46, p.96]

(1) $\qquad \varsigma(s) = \sum_{n=1}^{\infty}\frac{1}{n^s} = \frac{1}{1-2^{-s}}\sum_{n=1}^{\infty}\frac{1}{(2n-1)^s} = \frac{1}{1-2^{-s}}\sum_{n=0}^{\infty}\frac{1}{(2n+1)^s} \qquad ,(\operatorname{Re}(s) > 1)$

$$= \frac{1}{1-2^{1-s}} \sum_{n=1}^{\infty} \frac{(-1)^{n+1}}{n^s} = \frac{1}{1-2^{1-s}} \varsigma_a(s) \qquad , (\operatorname{Re}(s) > 0; \; s \neq 1)$$

where $\varsigma_a(s)$ is the alternating zeta function.

In their 1995 paper "Mellin Transforms and Asymptotics: Finite Differences and Rice's Integrals", Flajolet and Sedgewick [19] employed the following lemma:

Let $\varphi(z)$ be analytic in a domain that contains the half line $[n_0, \infty)$. Then the differences of the sequence $\{\varphi(k)\}$ admit the integral representation

(2) $$\sum_{k=0}^{n} \binom{n}{k} (-1)^k \varphi(k) = \frac{(-1)^n}{2\pi i} \int_C \varphi(z) \frac{n!}{z(z-1)\ldots(z-n)} dz$$

where $C$ is a positively oriented closed curve that lies in the domain of analyticity of $\varphi(z)$, encircles $[n_0, n)$ and does not include any of the integers $0, 1, \ldots, n_0 - 1$.

Using the above lemma with $\varphi(k) = \frac{1}{k^m}$, Flajolet and Sedgewick [19] also proved the following identity:

Defining $S_n(m)$ by

(3) $$S_n(m) = \sum_{k=1}^{n} \binom{n}{k} \frac{(-1)^k}{k^m}$$

then $S_n(m)$ can be expressed in terms of the generalised harmonic numbers as

(4) $$-S_n(m) = \sum_{1m_1 + 2m_2 + 3m_3 \ldots = m} \frac{1}{m_1! \, m_2! \, m_3! \ldots} \left(\frac{H_n^{(1)}}{1}\right)^{m_1} \left(\frac{H_n^{(2)}}{2}\right)^{m_2} \left(\frac{H_n^{(3)}}{3}\right)^{m_3} \ldots$$

where $H_n^{(k)}$ are the generalised harmonic numbers defined by

(5) $$H_n^{(k)} = \sum_{j=1}^{n} \frac{1}{j^k}$$

The first few values of $S_n(m)$ given by Flajolet and Sedgewick are:

(6.1) $$-S_n(1) = -\sum_{k=1}^{n} \binom{n}{k} \frac{(-1)^k}{k} = H_n^{(1)}$$

(6.2) $$-S_n(2) = -\sum_{k=1}^{n}\binom{n}{k}\frac{(-1)^k}{k^2} = \frac{1}{2}\left(H_n^{(1)}\right)^2 + \frac{1}{2}H_n^{(2)}$$

(6.3) $$-S_n(3) = -\sum_{k=1}^{n}\binom{n}{k}\frac{(-1)^k}{k^3} = \frac{1}{6}\left(H_n^{(1)}\right)^3 + \frac{1}{2}H_n^{(1)}H_n^{(2)} + \frac{1}{3}H_n^{(3)}$$

Adamchik's 1996 paper, "On Stirling Numbers and Euler Sums" [1] contained the following identities:

(7.1) $$\sum_{k=1}^{n}\frac{H_k^{(1)}}{k} = \frac{1}{2}\left(H_n^{(1)}\right)^2 + \frac{1}{2}H_n^{(2)}$$

(7.2) $$\sum_{k=1}^{n}\frac{H_k^{(2)}}{k} + \sum_{k=1}^{n}\frac{H_k^{(1)}}{k^2} = H_n^{(3)} + H_n^{(1)}H_n^{(2)}$$

(7.3) $$\sum_{k=1}^{n}\frac{\left(H_k^{(1)}\right)^2}{k} + \sum_{k=1}^{n}\frac{H_k^{(2)}}{k} = \frac{1}{3}\left(H_n^{(1)}\right)^3 + H_n^{(1)}H_n^{(2)} + \frac{2}{3}H_n^{(3)}$$

$$= 2\sum_{k=1}^{n}\frac{1}{k}\sum_{j=1}^{k}\frac{H_j^{(1)}}{j}$$

As a matter of interest, I also found formula (7.1) reported by Levenson in a 1938 volume of The American Mathematical Monthly [32] in a problem concerning the evaluation of

(8) $$\Gamma''(1) = \int_0^{\infty} e^{-x}\log^2 x\, dx = \gamma^2 + \varsigma(2)$$

where $\gamma$ is Euler's constant defined by

(9) $$\gamma = \lim_{n\to\infty}(H_n - \log n)$$

(and, regarding the above integral, refer to (24.2) below).

The third identity (7.3) is equal to $-2S_n(3)$ and hence we have

(10) $$S_n(3) = \sum_{k=1}^{n}\binom{n}{k}\frac{(-1)^k}{k^3}$$

$$= -\left\{\frac{1}{6}\left(H_n^{(1)}\right)^3 + \frac{1}{2}H_n^{(1)}H_n^{(2)} + \frac{1}{3}H_n^{(3)}\right\}$$

$$= -\frac{1}{2}\left\{\sum_{k=1}^{n}\frac{\left(H_k^{(1)}\right)^2}{k} + \sum_{k=1}^{n}\frac{H_k^{(2)}}{k}\right\}$$

The following formula was found by Sondow [43] by applying the Euler series transformation method (which is covered in Knopp's excellent book [28, p.240]) to the alternating Riemann zeta function. Sondow's result was

(11) $$\varsigma_a(s) = \sum_{n=0}^{\infty}\frac{1}{2^{n+1}}\sum_{k=0}^{n}\binom{n}{k}\frac{(-1)^k}{(k+1)^s}$$

where the alternating Riemann zeta function is defined by

$$\varsigma_a(s) = \sum_{n=1}^{\infty}\frac{(-1)^{n+1}}{n^s} = \sum_{n=0}^{\infty}\frac{(-1)^n}{(n+1)^s}$$

and is sometimes called the Dirichlet eta function and often designated by $\eta(s)$. It is known that $\varsigma_a(s)$ is an analytic function for $\text{Re}(s) > 0$.

The identity (11) has some history: it was conjectured by Knopp (1882-1957) around 1930, then proved by Hasse [25] in 1930 and subsequently rediscovered by Sondow in 1994. Hasse (1898-1979) also showed that

(12.1) $$\varsigma(s) = \frac{1}{s-1}\sum_{n=0}^{\infty}\frac{1}{n+1}\sum_{k=0}^{n}\binom{n}{k}\frac{(-1)^k}{(k+1)^{s-1}}$$

(12.2) $$\varsigma(s,x) = \frac{1}{s-1}\sum_{n=0}^{\infty}\frac{1}{n+1}\sum_{k=0}^{n}\binom{n}{k}\frac{(-1)^k}{(k+x)^{s-1}}$$

where $\varsigma(s,x) = \sum_{n=0}^{\infty}\frac{1}{(n+x)^s}$ for $\text{Re}(s) > 1$ is the Hurwitz zeta function. The above two formulae are valid for all $s$ except $s = 1$. It may be immediately noted that $\varsigma(s,1) = \varsigma(s)$.

It is shown in equation (4.4.85) in [14] that

(12.3) $$\frac{\log y}{s-1}Li_{s-1}(y) + Li_s(y) = \frac{y}{s-1}\sum_{n=0}^{\infty}\frac{1}{n+1}\sum_{k=0}^{n}\binom{n}{k}(-1)^k\frac{y^k}{(k+1)^{s-1}}$$

and with $y = 1$ this reverts to (12.1). A different proof of (12.2) has recently been given by Amore [3].

It is shown in equation (4.4.79) in [14] that

(13) $$\varsigma_a(s,x) = \sum_{n=0}^{\infty} \frac{1}{2^{n+1}} \sum_{k=0}^{n} \binom{n}{k} \frac{(-1)^k}{(k+x)^s}$$

where $\varsigma_a(s,u)$ may be regarded as an alternating Hurwitz zeta function and this may be written as

(13.1) $$\varsigma_a(s,x) = \sum_{n=0}^{\infty} \frac{(-1)^n}{(n+x)^s}$$

At first sight, the two Hasse identities (12.2) and (13) look rather different. However, if we consider the function defined by

$$f(t,s) = \sum_{n=0}^{\infty} t^n \sum_{k=0}^{n} \binom{n}{k} \frac{(-1)^k}{(k+1)^s}$$

we see that they are in fact intimately related. Indeed we have

$$(s-1)\varsigma(s) = \sum_{n=0}^{\infty} \frac{1}{n+1} \sum_{k=0}^{n} \binom{n}{k} \frac{(-1)^k}{(k+1)^{s-1}} = \int_0^1 f(t, s-1)\, dx$$

$$\varsigma_a(s) = \sum_{n=0}^{\infty} \frac{1}{2^{n+1}} \sum_{k=0}^{n} \binom{n}{k} \frac{(-1)^k}{(k+1)^s} = \frac{1}{2} f(1/2, s)$$

It is also shown in equation (4.4.99aiv) in [14] that

$$\sum_{n=1}^{\infty} t^n \sum_{k=0}^{n} \binom{n}{k} \frac{x^k}{(k+y)^s} = \frac{1}{\Gamma(s)} \int_0^{\infty} \frac{u^{s-1} e^{-yu} (1+xe^{-u}) t}{\left[1-(1+xe^{-u})t\right]}\, du$$

A number of identities were determined in [11] to [15] by using the Hasse formula, a combination of the Flajolet and Sedgewick formulae (6) and the Adamchik identities (7). Some of these identities are collected below for ease of reference.

$$\varsigma_a(s) = \sum_{n=0}^{\infty} \frac{1}{2^{n+1}} \sum_{k=0}^{n} \binom{n}{k} \frac{(-1)^k}{(k+1)^s}$$

$$\varsigma_a(1) = \sum_{n=1}^{\infty} \frac{1}{n 2^n} = \log 2$$

$$\varsigma_a(2) = \sum_{n=1}^{\infty} \frac{H_n}{n 2^n}$$

$$\varsigma_a(3) = \frac{1}{2} \sum_{n=1}^{\infty} \frac{1}{n 2^n} \left\{ (H_n^{(1)})^2 + H_n^{(2)} \right\} = \sum_{n=0}^{\infty} \frac{1}{n 2^n} \sum_{k=1}^{n} \binom{n}{k} \frac{(-1)^{k+1}}{k^2}$$

$$\varsigma_a(4) = \sum_{n=1}^{\infty} \frac{1}{n 2^n} \left\{ \frac{1}{6}(H_n^{(1)})^3 + \frac{1}{2} H_n^{(1)} H_n^{(2)} + \frac{1}{3} H_n^{(3)} \right\} = \sum_{n=1}^{\infty} \frac{1}{n 2^n} \sum_{k=1}^{n} \binom{n}{k} \frac{(-1)^{k+1}}{k^3}$$

$$\varsigma_a(5) = \sum_{n=1}^{\infty} \frac{1}{n 2^n} \left\{ (H_n^{(1)})^4 + 6(H_n^{(1)})^2 H_n^{(2)} + 8 H_n^{(1)} H_n^{(3)} + 3(H_n^{(2)})^2 + 6 H_n^{(4)} \right\}$$

$$= \sum_{n=1}^{\infty} \frac{1}{n 2^n} \sum_{k=1}^{n} \binom{n}{k} \frac{(-1)^{k+1}}{k^4}$$

$$\varsigma(s) = \frac{1}{s-1} \sum_{n=0}^{\infty} \frac{1}{n+1} \sum_{k=0}^{n} \binom{n}{k} \frac{(-1)^k}{(k+1)^{s-1}}$$

$$\varsigma(2) = \sum_{n=1}^{\infty} \frac{1}{n^2}$$

$$\varsigma(3) = \frac{1}{2!} \sum_{n=1}^{\infty} \frac{H_n}{n^2} = \frac{1}{2} \sum_{n=1}^{\infty} \frac{1}{n^2} \sum_{k=1}^{n} \binom{n}{k} \frac{(-1)^{k+1}}{k}$$

$$\varsigma(4) = \frac{1}{3!} \sum_{n=1}^{\infty} \frac{1}{n^2} \left\{ (H_n^{(1)})^2 + H_n^{(2)} \right\} = \frac{1}{3} \sum_{n=1}^{\infty} \frac{1}{n^2} \sum_{k=1}^{n} \binom{n}{k} \frac{(-1)^{k+1}}{k^2}$$

$$\varsigma(5) = \frac{1}{4!} \sum_{n=1}^{\infty} \frac{1}{n^2} \left\{ (H_n^{(1)})^3 + 3 H_n^{(1)} H_n^{(2)} + 2 H_n^{(3)} \right\} = \frac{1}{4} \sum_{n=1}^{\infty} \frac{1}{n^2} \sum_{k=1}^{n} \binom{n}{k} \frac{(-1)^{k+1}}{k^3}$$

On the basis of this limited data one may conjecture that

(14.1) $$\varsigma(s+2) = \frac{1}{s+1} \sum_{n=1}^{\infty} \frac{1}{n^2} \sum_{k=1}^{n} \binom{n}{k} \frac{(-1)^{k+1}}{k^s} \qquad , s \geq 1$$

(14.2) $$\varsigma_a(s) = \sum_{n=1}^{\infty} \frac{1}{n 2^n} \sum_{k=1}^{n} \binom{n}{k} \frac{(-1)^{k+1}}{k^{s-1}} \qquad , s \geq 1$$

These conjectures were proved (and in fact generalised) in equations (4.4.58) and (4.4.45) respectively of [14]: the generalised identities are

(14.3) $$\sum_{n=1}^{\infty}\frac{1}{n^2}\sum_{k=1}^{n}\binom{n}{k}\frac{(-1)^k x^k}{k^s} = -(s+1)Li_{s+2}(x) + \log x\, Li_{s+1}(x)$$

(14.4) $$\sum_{n=1}^{\infty}\frac{1}{n2^n}\sum_{k=1}^{n}\binom{n}{k}\frac{x^k}{k^s} = Li_{s+1}(x)$$

where the polylogarithm function $Li_n(x)$ is defined by

(14.5) $$Li_n(x) = \sum_{k=1}^{\infty}\frac{x^k}{k^n}\quad,\ (|x|\le 1)$$

We also note that Spieß [44] has derived the following identities

(15) $$\sum_{k=1}^{n}\frac{1}{k}\frac{1}{n-k+1} = \frac{2}{n+1}H_n^{(1)}$$

$$\sum_{k=2}^{n}\frac{2}{k}\frac{1}{n-k+1}H_{k-1}^{(1)} = \frac{3}{n+1}\left[\left(H_n^{(1)}\right)^2 - H_n^{(2)}\right]$$

$$\sum_{k=2}^{n}\frac{2}{k}\frac{2}{n-k+1}H_{k-1}^{(1)}H_{n-k}^{(1)} = \frac{4}{n+1}\left[\left(H_n^{(1)}\right)^3 - 3H_n^{(1)}H_n^{(2)} + 2H_n^{(3)}\right]$$

Reference should also be made to the paper by Larcombe et al. [31] where they show that for integers $m \ge 1, n \ge 0$.

(16.1) $$m\binom{m+n}{n}\sum_{k=0}^{n}\binom{n}{k}\frac{(-1)^k}{m+k} = 1$$

(16.2) $$m\binom{m+n}{n}\sum_{k=0}^{n}\binom{n}{k}\frac{(-1)^k}{(m+k)^2} = \sum_{k=m}^{m+n}\frac{1}{k}$$

(16.3) $$2m\binom{m+n}{n}\sum_{k=0}^{n}\binom{n}{k}\frac{(-1)^k}{(m+k)^3} = \left(\sum_{k=m}^{m+n}\frac{1}{k}\right)^2 + \sum_{k=m}^{m+n}\frac{1}{k^2}$$

(16.4) $$6m\binom{m+n}{n}\sum_{k=0}^{n}\binom{n}{k}\frac{(-1)^k}{(m+k)^4} = \left(\sum_{k=m}^{m+n}\frac{1}{k}\right)^3 + 3\left(\sum_{k=m}^{m+n}\frac{1}{k}\right)\left(\sum_{k=m}^{m+n}\frac{1}{k^2}\right) + 2\sum_{k=m}^{m+n}\frac{1}{k^3}$$

In their paper they employ integrals of the type

(17) $$\int_0^\infty x^p e^{-mx}(1-e^{-x})^n \, dx.$$

and, if we use the substitution $t = e^{-x}$, we can immediately see the relationship with equation (4.4.16) in [13]

(17.1) $$g^{(s-1)}(x) = \int_0^1 t^{x-1}(1-t)^n \log^s t \, dt$$

where $g(x)$ is defined by (and this function is also employed in (38.5) below)

(17.2) $$g(x) = \frac{n!}{x(1+x)\ldots(n+x)} = \frac{\Gamma(n+1)\Gamma(x)}{\Gamma(n+1+x)}$$

We also note the similarity with the following identity given by Anglesio in [5].

(17.3) $$\int_0^\infty \frac{e^{-ax}(1-e^{-x})^n}{x^r} \, dx = \frac{(-1)^r}{(r-1)!} \sum_{k=0}^n \binom{n}{k}(-1)^k (a+k)^{r-1} \log(a+k)$$

where $a \geq 0$ and $1 \leq r \leq n$ (except for $a = 0, r = 1$).

In fact, Anglesio's proof contains the identity for $a, p > 0$

(17.4) $$\int_0^\infty e^{-ax}(1-e^{-x})^n x^{p-1} dx = \Gamma(p) \sum_{k=0}^n \binom{n}{k} \frac{(-1)^k}{(a+k)^p}$$

Note that by letting $m = 1$ in (16.2) and (16.3), we obtain equations (4.4.127) and (4.4.130) respectively of [14]

$$\sum_{k=0}^n \binom{n}{k} \frac{(-1)^k}{(k+1)^2} = \frac{H_{n+1}}{n+1}$$

$$\sum_{k=0}^n \binom{n}{k} \frac{(-1)^k}{(k+1)^3} = \frac{1}{n+1} \sum_{k=1}^{n+1} \frac{H_k}{k} = \frac{1}{2(n+1)} \left\{ \left(H_{n+1}^{(1)}\right)^2 + H_{n+1}^{(2)} \right\}$$

In addition, Larcombe et al. also show how the above formulae can be derived using the identity from Gould's book, "Combinatorial Identities" [22]

(18) $$f(x+y) = y \binom{y+n}{n} \sum_{k=0}^n (-1)^k \binom{n}{k} \frac{f(x-k)}{y+k}$$

where $f(t)$ is a polynomial of degree $\leq n$. In [21], Gould attributes the formula (18) to Melzak [30]. Reference should also be made to the paper by Kirschenhofer [27]. Alternative proofs of the Larcombe at al. identities are given in [14]. Vermaseren [48] has also considered various harmonic sums and related integrals.

It is shown in [14] that the generalised harmonic numbers also feature in some known logarithmic integrals related to the derivative of the beta function

$$(19.1) \qquad -n\int_0^1 (1-t)^{n-1} \log t \, dt = H_n^{(1)}$$

$$(19.2) \qquad n\int_0^1 (1-t)^{n-1} \log^2 t \, dt = H_n^{(2)} + \left(H_n^{(1)}\right)^2$$

$$(19.3) \qquad -n\int_0^1 (1-t)^{n-1} \log^3 t \, dt = 6\left(\frac{1}{6}\left[H_n^{(1)}\right]^3 + \frac{1}{2}H_n^{(1)} H_n^{(2)} + \frac{1}{3}H_n^{(3)}\right)$$

and it is proved in equation (4.4.155zi) of [15] that

$$(19.4) \qquad (-1)^{p+1} n\int_0^1 (1-t)^{n-1} \log^p t \, dt = p!\sum_{k=1}^n \binom{n}{k}\frac{(-1)^k}{k^p}$$

Indeed, since $\int_0^1 (1-t)^{n-1} \log^k t \, dt = \int_0^1 t^{n-1} \log^k (1-t) \, dt$, one would automatically expect a connection with the Stirling numbers of the first kind defined below in (20.2).

After that lengthy introduction, we now need to recall some properties of the Stirling numbers $s(n,k)$ of the first kind and the (exponential) complete Bell polynomials. These are considered in the next two sections.

## STIRLING NUMBERS OF THE FIRST KIND

The Stirling numbers $s(n,k)$ of the first kind [46, p.56] are defined by the generating function

$$(20.1) \qquad x(x-1)\ldots(x-n+1) = \sum_{k=0}^n s(n,k)x^k$$

and also by the Maclaurin expansion due to Cauchy

(20.2) $$\log^k(1+x) = k! \sum_{n=k}^{\infty} s(n,k) \frac{x^n}{n!} \qquad |x| < 1$$

Since $s(n,k) = \frac{1}{k!} \frac{d^n}{dx^n} \log^k(1+x) \Big|_{x=0}$ it is clear that $s(n,k) = 0 \ \forall \ n \leq k-1$ (as is also evident from the polynomial expression in (20.1)). We also have $s(0,0) = 1$.

The following proof of (20.2) was given by Póyla and Szegö in [37, p.227]. From the binomial theorem we have

$$(1-t)^{-x} = \sum_{n=0}^{\infty} \frac{x(x+1)(x+2)\ldots(x+n-1)}{n!} t^n$$

Using (20.1) this becomes

$$= \sum_{n=0}^{\infty} \frac{(-1)^n t^n}{n!} \sum_{k=0}^{n} s(n,k)(-1)^k x^k$$

$$= 1 + \sum_{n=1}^{\infty} \frac{(-1)^n t^n}{n!} \sum_{k=0}^{n} s(n,k)(-1)^k x^k$$

$$= 1 + \sum_{n=1}^{\infty} \frac{(-1)^n t^n}{n!} \sum_{k=1}^{n} s(n,k)(-1)^k x^k$$

$$= 1 + \sum_{k=1}^{\infty} (-1)^k x^k \sum_{n=k}^{\infty} \frac{s(n,k)}{n!} (-1)^n t^n$$

We then note that

$$(1-t)^{-x} = \exp[-x \log(1-t)]$$

$$= \sum_{k=1}^{\infty} \frac{x^k}{k!} \left( \log \frac{1}{1-t} \right)^k$$

and upon comparing coefficients of $x^k$ we obtain (20.2).

The first few Stirling numbers $s(n,k)$ of the first kind are given in [41] and also in the book by Srivastava and Choi [46, p.57]

(20.3) $\qquad s(n,0) = \delta_{n,0}$

$$s(n,1) = (-1)^{n+1}(n-1)!$$

$$s(n,2) = (-1)^n (n-1)! H_{n-1}$$

$$s(n,3) = (-1)^{n+1} \frac{(n-1)!}{2} \left\{ (H_{n-1})^2 - H_{n-1}^{(2)} \right\}$$

$$s(n,4) = (-1)^n \frac{(n-1)!}{6} \left\{ (H_{n-1})^3 - 3H_{n-1} H_{n-1}^{(2)} + 2H_{n-1}^{(3)} \right\}$$

The above representations should be compared with the identities found by Larcombe et al. as set out above in (16). An elementary method for determining the Stirling numbers $s(n,k)$ of the first kind is shown below in (44.6) et seq.

## COMPLETE BELL POLYNOMIALS

The (exponential) complete Bell polynomials may be defined by $Y_0 = 1$ and for $n \geq 1$

(21) $$Y_n(x_1,\ldots,x_n) = \sum_{\pi(n)} \frac{n!}{k_1! k_2! \ldots k_n!} \left(\frac{x_1}{1!}\right)^{k_1} \left(\frac{x_2}{2!}\right)^{k_2} \ldots \left(\frac{x_n}{n!}\right)^{k_n}$$

where the sum is taken over all partitions $\pi(n)$ of $n$, i.e. over all sets of integers $k_j$ such that

(21.1) $$k_1 + 2k_2 + 3k_3 + \ldots + nk_n = n$$

The complete Bell polynomials have integer coefficients and the first six are set out below (Comtet [10, p.307])

(22) $$Y_1(x_1) = x_1$$

$$Y_2(x_1, x_2) = x_1^2 + x_2$$

$$Y_3(x_1, x_2, x_3) = x_1^3 + 3x_1 x_2 + x_3$$

$$Y_4(x_1, x_2, x_3, x_4) = x_1^4 + 6x_1^2 x_2 + 4x_1 x_3 + 3x_2^2 + x_4$$

$$Y_5(x_1, x_2, x_3, x_4, x_5) = x_1^5 + 10x_1^3 x_2 + 10x_1^2 x_3 + 15x_1 x_2^2 + 5x_1 x_4 + 10x_2 x_3 + x_5$$

$$Y_6(x_1, x_2, x_3, x_4, x_5, x_6) = x_1^6 + 6x_1 x_5 + 15x_2 x_4 + 10x_2^3 + 15x_1^2 x_4 + 15x_2^3 + 60x_1 x_2 x_3$$

$$+ 20x_1^3 x_3 + 45x_1^2 x_2^2 + 15x_1^4 x_1 + x_6$$

The total number of terms $\pi(n)$ increases rapidly; for example, as reported by Bell [6] in 1934, we have $\pi(22) = 1002$ terms.

The modus operandi of the summation in (21.1) is easily illustrated by the following example: if $n = 4$, then $k_1 + 2k_2 + 3k_3 + 4k_4 = 4$ is satisfied by the integers in the following array

$$\begin{bmatrix} k_1 & k_2 & k_3 & k_4 \\ 4 & 0 & 0 & 0 \\ 2 & 1 & 0 & 0 \\ 1 & 0 & 1 & 0 \\ 0 & 2 & 0 & 0 \\ 0 & 0 & 0 & 1 \end{bmatrix}$$

and these powers feature in the series below

$$\varsigma_a(5) = \sum_{n=1}^{\infty} \frac{1}{n 2^n} \left\{ \left(H_n^{(1)}\right)^4 + 6\left(H_n^{(1)}\right)^2 H_n^{(2)} + 8 H_n^{(1)} H_n^{(3)} + 3\left(H_n^{(2)}\right)^2 + 6 H_n^{(4)} \right\}$$

The complete Bell polynomials are also given by the exponential generating function (Comtet [10, p.134])

(23) $$\exp\left(\sum_{j=1}^{\infty} x_j \frac{t^j}{j!}\right) = 1 + \sum_{n=1}^{\infty} Y_n(x_1, \ldots, x_n) \frac{t^n}{n!} = \sum_{n=0}^{\infty} Y_n(x_1, \ldots, x_n) \frac{t^n}{n!}$$

Let us now consider a function $f(x)$ which has a Taylor series expansion around $x$: we have

$$e^{f(x+t)} = \exp\left(\sum_{j=0}^{\infty} f^{(j)}(x) \frac{t^j}{j!}\right) = e^{f(x)} \exp\left(\sum_{j=1}^{\infty} f^{(j)}(x) \frac{t^j}{j!}\right)$$

$$= e^{f(x)} \left\{ 1 + \sum_{n=1}^{\infty} Y_n\left(f^{(1)}(x), f^{(2)}(x), \ldots, f^{(n)}(x)\right) \frac{t^n}{n!} \right\}$$

We see that

$$\frac{d^m}{dx^m} e^{f(x)} = \frac{\partial^m}{\partial x^m} e^{f(x+t)} \bigg|_{t=0} = \frac{\partial^m}{\partial t^m} e^{f(x+t)} \bigg|_{t=0}$$

and we therefore obtain (as noted by Kölbig [30] and Coffey [8])

(24) $$\frac{d^m}{dx^m} e^{f(x)} = e^{f(x)} Y_m\left(f^{(1)}(x), f^{(2)}(x),..., f^{(m)}(x)\right)$$

Differentiating (24) we see that

$$Y_{m+1}\left(f^{(1)}(x), f^{(2)}(x),..., f^{(m+1)}(x)\right) = \left(f^{(1)}(x) + \frac{d}{dx}\right) Y_m\left(f^{(1)}(x), f^{(2)}(x),..., f^{(m)}(x)\right)$$

As an example, letting $f(x) = \log \Gamma(x)$ in (24) we obtain

(24.1) $$\frac{d^m}{dx^m} e^{\log \Gamma(x)} = \Gamma^{(m)}(x) = \Gamma(x) Y_m\left(\psi(x), \psi^{(1)}(x),..., \psi^{(m-1)}(x)\right)$$

$$= \int_0^\infty t^{x-1} e^{-t} \log^m t \, dt$$

and since [46, p.22]

$$\psi^{(p)}(x) = (-1)^{p+1} p! \varsigma(p+1, x)$$

we may express $\Gamma^{(m)}(x)$ in terms of $\psi(x)$ and the Hurwitz zeta functions. In particular, Coffey [8] notes that

(24.2) $$\Gamma^{(m)}(1) = Y_m(-\gamma, x_1,..., x_{m-1})$$

where $x_p = (-1)^{p+1} p! \varsigma(p+1)$. Values of $\Gamma^{(m)}(1)$ are reported in [46, p.265] for $m \leq 10$ and the first three are

$$\Gamma^{(1)}(1) = -\gamma$$

$$\Gamma^{(2)}(1) = \varsigma(2) + \gamma^2$$

$$\Gamma^{(3)}(1) = -2\varsigma(3) - 3\gamma\varsigma(2) - \gamma^3$$

We have [46, p.20]

$$\psi\left(\frac{1}{2}\right) = -\gamma - 2\log 2$$

$$\psi^{(m)}\left(\frac{1}{2}\right) = (-1)^{m+1} m!(2^{m+1} - 1)\varsigma(m+1)$$

and therefore we may readily obtain an expression for $\Gamma^{(m)}\left(\frac{1}{2}\right)$. The first ten values are also reported in [46, p.266]. We could also, for example, let $f(x) = \log\sin(\pi x)$ in (24) to obtain complete Bell polynomials involving the derivatives of $\cot(\pi x)$.

The following is extracted from a series of papers written by Snowden [42, p.68] in 2003.

Let us consider the function $f(x)$ with the following Maclaurin expansion

(24.3) $$\log f(x) = b_0 + \sum_{n=1}^{\infty} \frac{b_n}{n} x^n$$

and we wish to determine the coefficients $a_n$ such that

$$f(x) = \sum_{n=0}^{\infty} a_n x^n$$

By differentiating (24.3) and multiplying the two power series, we get

$$na_n = \sum_{k=1}^{n} b_k a_{n-k}$$

Upon examination of this recurrence relation it is easy to see that

$$n!a_n = a_0[b_1, -b_2, b_3, ..., (-1)^{n+1} b_n]$$

where the symbol $[a_1, a_2, a_3, ..., a_n]$ is defined as the $n \times n$ determinant

$$\begin{vmatrix} a_1 & a_2 & a_3 & a_4 & . & . & . & a_n \\ (n-1) & a_1 & a_2 & a_3 & . & . & . & a_{n-1} \\ 0 & (n-2) & a_1 & a_2 & . & . & . & a_{n-2} \\ 0 & 0 & (n-3) & a_1 & . & . & . & a_{n-3} \\ \vdots & \vdots & \vdots & \vdots & \vdots & \vdots & \vdots & \vdots \\ 0 & 0 & 0 & 0 & 0 & 0 & 1 & a_1 \end{vmatrix}$$

Since $\log f(0) = \log a_0 = b_0$ we have

$$f(x) = e^{b_0}\left[1 + \sum_{n=1}^{\infty} [b_1, -b_2, b_3, ..., (-1)^{n+1} b_n] \frac{x^n}{n!}\right]$$

Multiplying (24.3) by $\alpha$ it is easily seen that

(24.4) $$f^{\alpha}(x) = e^{-b_0}\left[1 + \sum_{n=1}^{\infty}[\alpha b_1, -\alpha b_2, \alpha b_3, ..., (-1)^{n+1}\alpha b_n]\frac{x^n}{n!}\right]$$

and, in particular, with $\alpha = -1$ we obtain

(24.5) $$\frac{1}{f(x)} = e^{-b_0}\left[1 + \sum_{n=1}^{\infty}[-b_1, b_2, -b_3, ..., (-1)^n b_n]\frac{x^n}{n!}\right]$$

Differentiating (24.4) with respect to $\alpha$ would give us an expression for $f^{\alpha}(x)\log f(x)$.

We note the well known series expansion (which is also derived as a by-product later in equation (46))

$$\log \Gamma(1+x) = -\gamma x + \sum_{n=2}^{\infty}(-1)^n \frac{\varsigma(n)}{n} x^n \quad , -1 < x \leq 1$$

and hence we have

(24.6) $$\Gamma(1+x) = 1 + \sum_{n=1}^{\infty}[-\gamma, -\varsigma(2), -\varsigma(3), ..., -\varsigma(n)]\frac{x^n}{n!}$$

and

(24.7) $$\frac{1}{\Gamma(1+x)} = 1 + \sum_{n=1}^{\infty}[\gamma, \varsigma(2), \varsigma(3), ..., \varsigma(n)]\frac{x^n}{n!}$$

where $\varsigma(1)$ is defined as equal to $\gamma$.

Differentiating (24.6) we get

$$\Gamma^{(m)}(1+x) = \sum_{n=1}^{\infty}[-\gamma, -\varsigma(2), -\varsigma(3), ..., -\varsigma(n)]\frac{n(n-1)...(n-m+1)x^{n-m}}{n!}$$

and hence, letting $x = 0$, we have the $m$ th derivative of the gamma function in the form of a determinant

(24.8) $$\Gamma^{(m)}(1) = [-\varsigma(1), -\varsigma(2), -\varsigma(3), ..., -\varsigma(m)]$$

where we again designate $\varsigma(1) = \gamma$.

In his paper "The asymptotic behaviour of the Stirling numbers of the first kind" [50], Wilf proved that if $\begin{bmatrix} n \\ k \end{bmatrix}$ is the (signless) Stirling number of the first kind, then for each fixed integer $k \geq 2$ we have (the signless Stirling number of the first kind is defined as the absolute value of $s(n,k)$)

$$(24.9) \quad \frac{1}{(n-1)!}\begin{bmatrix} n \\ k \end{bmatrix} = \lambda_1 \frac{\log^{k-1} n}{(k-1)!} + \lambda_2 \frac{\log^{k-2} n}{(k-2)!} + \ldots + \lambda_k + O\left(\frac{\log^{k-2} n}{n}\right)$$

where $\lambda_j$ are the coefficients in the expansion

$$(24.10) \quad \frac{1}{\Gamma(x)} = \sum_{j=1}^{\infty} \lambda_j x^j$$

Since $\frac{1}{\Gamma(1+x)} = \frac{1}{x\Gamma(x)}$ equation (24.7) may be compared with (24.10)

In particular, using the recurrence (obtained from the logarithmic derivative of the infinite product for $\Gamma(z)$)

$$(24.11) \quad \lambda_{n+1} = \frac{1}{n}\left\{\gamma \lambda_n + \sum_{j=0}^{n-2}(-1)^{n-j-1}\varsigma(n-j)\lambda_{j+1}\right\} \quad n \geq 1; \lambda_1 = 1$$

we obtain

$$\lambda_1 = 1$$

$$\lambda_2 = \gamma$$

$$\lambda_3 = \frac{1}{12}(6\gamma^2 - \pi^2)$$

$$\lambda_4 = \frac{1}{12}\left[2\gamma^3 - \gamma\pi^2 + 4\varsigma(3)\right]$$

$$\lambda_5 = \frac{1}{1440}\left[60\gamma^4 - 60\gamma^2\pi^2 + \pi^4 + 480\varsigma(3)\right]$$

This material is considered further in [17].

Often in mathematics we look for divine inspiration but we do not usually expect to obtain it from a canonised saint. This is indeed the source of the next remark. The calculation of $g^{(m)}(x)$, as defined above in (17.2), effectively involves the derivative of a composite function $g(f(t))$ and the general formula for this was discovered by Francesco Faà di Bruno (1825-1888) who was declared a Saint by Pope John Paul II in St. Peter's Square in Rome in 1988 [43].

In [26] di Bruno showed that

(25) $$\frac{d^{(m)}}{dt^{(m)}} g(f(t)) = \sum \frac{m!}{b_1! b_2! \ldots b_m!} g^{(k)}(f(t)) \left(\frac{f^{(1)}(t)}{1!}\right)^{b_1} \left(\frac{f^{(2)}(t)}{2!}\right)^{b_2} \ldots \left(\frac{f^{(m)}(t)}{m!}\right)^{b_m}$$

where the sum is over all different solutions in non-negative integers $b_1, \ldots, b_m$ of $b_1 + 2b_2 + \ldots + mb_m$, and $k = b_1 + \ldots + b_m$. In our case, the composite function was of the form $g(f(t))$ where $g(t) = 1/t$ and $f(t) = t(t+1)\ldots(t+n)$ and we therefore have an explanation why the Flajolet and Sedgewick formulation of (4) and the di Bruno formula both involve the use of the complete Bell polynomials.

In [45] Gould reminded the mathematical community of the "not well-known" formula for the $n$ th derivative of a composite function $f(z)$ where $z$ is a function of $x$, namely

(26.1) $$D_x^{(n)} f(z) = \sum_{k=1}^{n} D_z^{(k)} f(z) \frac{(-1)^k}{k!} \sum_{j=1}^{k} (-1)^j \binom{k}{j} z^{k-j} \quad , \text{ for } n \geq 1$$

This expression is frequently easier to handle than the di Bruno algorithm.

This differentiation algorithm was also reported by Gould [23] in 1972 in a somewhat different form in the case where $z$ is a function of $s$

(26.2) $$D_s^{(n)} f(z) = \sum_{k=0}^{n} D_z^{(k)} f(z) \frac{(-1)^k}{k!} \sum_{j=0}^{k} (-1)^j \binom{k}{j} z^{k-j} D_s^{(n)} z^j \quad , \text{ for } n \geq 1$$

An interesting paper on the role of Bell polynomials in integration has recently been presented by Collins [9].

We are now ready to expand upon the work originally carried out by Coppo [18].

COPPO'S FORMULA

The following proposition was derived by Coppo [18] in 2003. For positive integers $q \geq 1$ and complex values of $x \in \mathbf{C}$, $x \neq 0, -1, \ldots, -n$ we have

(27) $$\sum_{k=0}^{n}\binom{n}{k}\frac{(-1)^k}{(k+x)^q}=\frac{n!}{x(x+1)...(x+n)}\vartheta(n+1,q-1,x)$$

where

(27.1) $$\vartheta(n,m,x)=\frac{1}{m!}Y_m\left(0!H_n^{(1)}(x),1!H_n^{(2)}(x),...,(m-1)!H_n^{(m)}(x)\right)$$

and where $H_n^{(m)}(x)$ is the generalised harmonic number function defined by

(28) $$H_n^{(m)}(x)=\sum_{k=0}^{n-1}\frac{1}{(k+x)^m}$$

Flajolet and Sedgewick [19] refer to $H_n^{(m)}(x)$ as the incomplete Hurwitz zeta functions (and hence the title of this paper). We note that

$$H_n^{(m)}(1)=H_n^{(m)}=\sum_{k=1}^{n}\frac{1}{k^m}$$

where $H_n^{(m)}$ are the generalised harmonic numbers.

This may be proved as follows:

With the definitions

$$S(n,0,x)=\frac{n!}{x(x+1)...(x+n-1)}$$

and

$$S(n,m,x)=\frac{(-1)^m}{m!}\frac{d^m}{dx^m}S(n,0,x)$$

we form the Taylor series

(29.1) $$S(n,0,x-t)=\sum_{m=0}^{\infty}S(n,m,x)t^m$$

We may write

$$S(n,0,x-t)=\frac{n!}{(x-t)(x-t+1)...(x-t+n-1)}$$

$$= \frac{n! \dfrac{1}{x} \dfrac{1}{1+x} \cdots \dfrac{1}{n-1+x}}{(1-t/x)(1-t/(1+x))\cdots(1-t/(n-1+x))}$$

$$= \frac{S(n,0,x)}{(1-t/x)(1-t/(1+x))\cdots(1-t/(n-1+x))}$$

and we may express this as

$$S(n,0,x-t) = S(n,0,x)\exp\left[-\sum_{k=0}^{n-1}\log(1-t/(k+x))\right]$$

Then, employing the Maclaurin expansion for $\log(1-t/(k+x))$, we obtain

$$S(n,0,x-t) = S(n,0,x)\exp\left[\sum_{k=0}^{n-1}\sum_{m=1}^{\infty}\frac{t^m}{m(k+x)}\right]$$

and reversing the order of summation this becomes

$$S(n,0,x-t) = S(n,0,x)\exp\left[\sum_{m=1}^{\infty}H_n^{(m)}(x)\frac{t^m}{m}\right]$$

(29.2) $\qquad = S(n,0,x)\exp\left[\sum_{m=1}^{\infty}(m-1)!\,H_n^{(m)}(x)\frac{t^m}{m!}\right]$

We now recall the definition of the complete Bell polynomials (23)

$$\exp\left[\sum_{m=1}^{\infty}x_m\frac{t^m}{m!}\right] = \sum_{m=0}^{\infty}Y_m(x_1,\ldots,x_m)\frac{t^m}{m!}$$

and see that

$$\exp\left[\sum_{m=1}^{\infty}H_n^{(m)}(x)\frac{t^m}{m}\right] = \sum_{m=0}^{\infty}Y_m\left(0!H_n^{(1)}(x),\,1!H_n^{(2)}(x),\ldots,(m-1)!H_n^{(m)}(x)\right)\frac{t^m}{m!}$$

Using (29.1) above

$$S(n,0,x-t) = \sum_{m=0}^{\infty}S(n,m,x)t^m$$

and equating coefficients of $t^m$ we see that

(29.3) $\qquad S(n,m,x) = S(n,0,x)\dfrac{1}{m!} Y_m\left(0!H_n^{(1)}(x), 1!H_n^{(2)}(x),...,(m-1)!H_n^{(m)}(x)\right)$

We have the well-known formula (see for example [12] and [24])

$$\dfrac{n!}{x(x+1)...(x+n)} = \sum_{k=0}^{n} \binom{n}{k} \dfrac{(-1)^k}{k+x}$$

which may be written as

$$\dfrac{1}{n+1} S(n+1,0,x) = \sum_{k=0}^{n} \binom{n}{k} \dfrac{(-1)^k}{k+x}$$

and differentiating this $q-1$ times gives us

(29.4) $\qquad \dfrac{1}{n+1} S(n+1, q-1, x) = \sum_{k=0}^{n} \binom{n}{k} \dfrac{(-1)^k}{(k+x)^q}$

Then, equating (29.3) and (29.4) gives us (27).

Since

$$\dfrac{n!}{x(x+1)...(x+n)} = \dfrac{\Gamma(n+1)\Gamma(x)}{\Gamma(n+1+x)}$$

we may write Coppo's formula as

(30)

$$\sum_{k=0}^{n} \binom{n}{k} \dfrac{(-1)^k}{(k+x)^q} = \dfrac{\Gamma(n+1)\Gamma(x)}{\Gamma(n+1+x)} \dfrac{1}{(q-1)!} Y_{q-1}\left(0!H_{n+1}^{(1)}(x), 1!H_{n+1}^{(2)}(x), ..., (q-2)!H_{n+1}^{(q-1)}(x)\right)$$

Designating $h(x) = H_{n+1}^{(1)}(x)$ we see that $h^{(p)}(x) = (-1)^p\, p!\, H_{n+1}^{(p+1)}(x)$ and we can therefore write (30) as

(30.1) $\qquad \sum_{k=0}^{n} \binom{n}{k} \dfrac{(-1)^k}{(k+x)^q} = \dfrac{\Gamma(n+1)\Gamma(x)}{\Gamma(n+1+x)} \dfrac{1}{(q-1)!} Y_{q-1}\left(h(x), -h^{(1)}(x), ..., (-1)^q h^{(q)}(x)\right)$

Coffey [8] also reported a version of this formula in 2006.

# AN APPLICATION OF THE HASSE IDENTITY FOR THE HURWITZ ZETA FUNCTION

As mentioned above, Hasse (1898-1979) showed that

$$(31.1) \qquad \varsigma(s) = \frac{1}{s-1} \sum_{n=0}^{\infty} \frac{1}{n+1} \sum_{k=0}^{n} \binom{n}{k} \frac{(-1)^k}{(k+1)^{s-1}}$$

$$(31.2) \qquad \varsigma(s,x) = \frac{1}{s-1} \sum_{n=0}^{\infty} \frac{1}{n+1} \sum_{k=0}^{n} \binom{n}{k} \frac{(-1)^k}{(k+x)^{s-1}}$$

where $\varsigma(s,x) = \sum_{n=0}^{\infty} \frac{1}{(n+x)^s}$ for $\operatorname{Re}(s) > 1$ is the Hurwitz zeta function. The above two formulae are valid for all $s$ except $s = 1$. It may be noted that $\varsigma(s,1) = \varsigma(s)$.

Using (30) and the Hasse identity (31.2) we obtain for $s = q+1$ and $q \geq 1$

$$(32)$$
$$\varsigma(q+1,x) = \frac{1}{q!} \sum_{n=0}^{\infty} \frac{1}{n+1} \frac{\Gamma(n+1)\Gamma(x)}{\Gamma(n+1+x)} Y_{q-1}\left(0! H_{n+1}^{(1)}(x), 1! H_{n+1}^{(2)}(x), \ldots, (q-2)! H_{n+1}^{(q-1)}(x)\right)$$

$$= \frac{1}{q!} \sum_{n=1}^{\infty} \frac{1}{n} \frac{\Gamma(n)\Gamma(x)}{\Gamma(n+x)} Y_{q-1}\left(0! H_{n}^{(1)}(x), 1! H_{n}^{(2)}(x), \ldots, (q-2)! H_{n}^{(q-1)}(x)\right)$$

Particular cases of Coppo's formula are set out below.

$$(33.1) \qquad \sum_{k=0}^{n} \binom{n}{k} \frac{(-1)^k}{(k+x)} = \frac{\Gamma(n+1)\Gamma(x)}{\Gamma(n+1+x)} Y_0 = \frac{\Gamma(n+1)\Gamma(x)}{\Gamma(n+1+x)}$$

$$(33.2) \qquad \sum_{k=0}^{n} \binom{n}{k} \frac{(-1)^k}{(k+x)^2} = \frac{\Gamma(n+1)\Gamma(x)}{\Gamma(n+1+x)} Y_1\left(0! H_{n+1}^{(1)}(x)\right) = \frac{\Gamma(n+1)\Gamma(x)}{\Gamma(n+1+x)} H_{n+1}^{(1)}(x)$$

$$\sum_{k=0}^{n} \binom{n}{k} \frac{(-1)^k}{(k+x)^3} = \frac{\Gamma(n+1)\Gamma(x)}{\Gamma(n+1+x)} \frac{1}{2} Y_2\left(0! H_{n+1}^{(1)}(x), 1! H_{n+1}^{(2)}(x)\right)$$

$$(33.3) \qquad = \frac{\Gamma(n+1)\Gamma(x)}{\Gamma(n+1+x)} \frac{1}{2}\left(\left[H_{n+1}^{(1)}(x)\right]^2 + H_{n+1}^{(2)}(x)\right)$$

$$\sum_{k=0}^{n} \binom{n}{k} \frac{(-1)^k}{(k+x)^4} = \frac{\Gamma(n+1)\Gamma(x)}{\Gamma(n+1+x)} \frac{1}{6} Y_3\left(H_{n+1}^{(1)}(x), H_{n+1}^{(2)}(x), 2 H_{n+1}^{(3)}(x)\right)$$

(33.4) $$= \frac{\Gamma(n+1)\Gamma(x)}{\Gamma(n+1+x)} \frac{1}{6}\left(\left[H_{n+1}^{(1)}(x)\right]^3 + 3H_{n+1}^{(1)}(x)H_{n+1}^{(2)}(x) + 2H_{n+1}^{(3)}(x)\right)$$

We have for $q=1$

(34) $$\varsigma(2,x) = \sum_{n=0}^{\infty} \frac{1}{n+1} \frac{\Gamma(n+1)\Gamma(x)}{\Gamma(n+1+x)} = \sum_{n=1}^{\infty} \frac{1}{n} \frac{\Gamma(n)\Gamma(x)}{\Gamma(n+x)}$$

and with $x=1$ this simply becomes the Riemann zeta function $\varsigma(2)$

$$\varsigma(2,1) = \varsigma(2) = \sum_{n=0}^{\infty} \frac{1}{(n+1)^2} = \sum_{n=1}^{\infty} \frac{1}{n^2}$$

Since $\varsigma(2,x) = \psi'(x)$ we also have

(35) $$\psi'(x) = \sum_{n=0}^{\infty} \frac{1}{n+1} \frac{\Gamma(n+1)\Gamma(x)}{\Gamma(n+1+x)}$$

which is reminiscent of the well-known formula

(36) $$\psi(x+a) - \psi(a) = \sum_{k=1}^{\infty} \frac{(-1)^{k+1}}{k} \frac{x(x-1)\dots(x-k+1)}{a(a+1)\dots(a+k-1)}$$

which converges for $\mathrm{Re}(x+a) > 0$. According to Raina and Ladda [38], this summation formula is due to Nörlund (see [35], [36]] and also Ruben's note [40]).

In equation (4.3.35) in [12] it was shown that

$$\psi'(x) = \sum_{n=1}^{\infty} \frac{1}{k} \frac{(n-1)!}{x(x+1)\dots(x+n-1)}$$

(36.1) $$= \sum_{n=1}^{\infty} \frac{\Gamma(n)}{n} \frac{\Gamma(x)}{\Gamma(x+n)}$$

This result was also reported by Ruben [40] in 1976. We see that (36.1) is equivalent to equation (34). This correspondence is considered further in (45) below.

We have Legendre's duplication formula [46, p.7] for $t > 0$

(36.2) $$\Gamma(t)\Gamma\left(t + \frac{1}{2}\right) = \frac{\sqrt{\pi}}{2^{2t-1}} \Gamma(2t)$$

Letting $t = n+1$ we obtain

$$\Gamma\left(n+1+\frac{1}{2}\right) = \frac{\sqrt{\pi}}{2^{2n+1}} \frac{\Gamma(2n+2)}{\Gamma(n+1)} = \frac{\sqrt{\pi}}{2^{2n+1}} \frac{(2n+1)!}{n!}$$

and we see that

$$\frac{\Gamma(n+1)\Gamma\left(\frac{1}{2}\right)}{\Gamma\left(n+1+\frac{1}{2}\right)} = \frac{2^{2n+1}(n!)^2}{(2n+1)!}$$

Hence, letting $x = 1/2$ in (34) we get

$$\varsigma\left(2,\frac{1}{2}\right) = \sum_{n=0}^{\infty} \frac{1}{n+1} \frac{2^{2n+1}(n!)^2}{(2n+1)!}$$

It is well known that

$$\varsigma\left(s,\frac{1}{2}\right) = (2^s - 1)\varsigma(s)$$

and hence we obtain

(37) $$\varsigma(2) = \frac{1}{3}\sum_{n=0}^{\infty} \frac{1}{n+1} \frac{2^{2n+1}[n!]^2}{(2n+1)!}$$

This has some structural similarities to Ramanujan's formula for Catalan's constant $G$ (see equation (8.15a) in [16], Adamchik's paper [2] and [36])

(38) $$G = \frac{\pi}{4}\sum_{n=0}^{\infty} \frac{1}{(2n+1)2^{4n}}\binom{2n}{n}^2 = \frac{\pi}{4}\sum_{n=0}^{\infty} \frac{[(2n)!]^2}{(2n+1)[(n)!]^4 \, 2^{4n}}$$

but the connection is more obvious with the expression contained in Ramanujan's Notebooks (Berndt [6i, Part I, p.289])

(38.1) $$G = \frac{1}{2}\sum_{n=0}^{\infty} \frac{1}{(2n+1)^2} \frac{2^{2n}[n!]^2}{(2n)!} = \frac{1}{4}\sum_{n=0}^{\infty} \frac{1}{(2n+1)} \frac{2^{2n+1}[n!]^2}{(2n+1)!}$$

In a personal communication Coffey mentioned that it seems that equations (38) and (38.1) have to be equivalent due to the duplication formula (36.2) for the gamma function.

About four years ago, the author derived the following expression for $G$ (see equation (6.30x) of [16])

$$G = \frac{1}{2}\int_0^1 \frac{\tan^{-1} x}{x}\,dx = \frac{1}{2}\sum_{n=0}^{\infty} \frac{2^{2n}(n!)^2}{(2n+1)!}\sum_{k=0}^{2n}(-1)^k \binom{2n}{k}\frac{2^{n-k}-1}{n-k}$$

and it is clear that this is equivalent to (38.1) provided one can prove that

$$\frac{1}{2n+1} = \sum_{k=0}^{2n}(-1)^k \binom{2n}{k}\frac{2^{n-k}-1}{n-k}$$

Using (37) and (38.1) we may write

$$3\varsigma(2) - 4G = \sum_{n=0}^{\infty}\frac{n}{(n+1)(2n+1)}\frac{2^{2n+1}[n!]^2}{(2n+1)!} = \sum_{n=0}^{\infty}\frac{n 2^{2n+1}}{(n+1)(2n+1)^2}\binom{2n}{n}^{-1}$$

which may possibly be of assistance in determining whether the difference $3\varsigma(2) - 4G$ is rational or irrational.

We have from (34)

$$\varsigma(2,x) = \sum_{n=1}^{\infty}\frac{\Gamma(n)}{n}\frac{\Gamma(x)}{\Gamma(n+x)}$$

$$= \sum_{n=1}^{\infty}\frac{B(n,x)}{n} = \sum_{n=1}^{\infty}\frac{1}{n}\int_0^1 t^{n-1}(1-t)^{x-1}\,dt$$

and we therefore obtain

(38.2) $$\varsigma(2,x) = -\int_0^1 \frac{(1-t)^{x-1}\log(1-t)}{t}\,dt$$

In particular we see that

$$\varsigma\left(2,\frac{1}{4}\right) = -\int_0^1 \frac{\log(1-t)}{t(1-t)^{3/4}}\,dt$$

In a scintilla temporis, the Wolfram Integrator miraculously tells us that

$$\frac{1}{4}\int \frac{\log(1-t)}{t(1-t)^{3/4}} dt = -Li_2\left(1-\sqrt[4]{1-t}\right) - \log\left(1+\sqrt[4]{1-t}\right)\log\left(\sqrt[4]{1-t}\right) - Li_2\left(-\sqrt[4]{1-t}\right)$$

$$+i\log\left(1+i\sqrt[4]{1-t}\right)\log\left(\sqrt[4]{1-t}\right) + iLi_2\left(-i\sqrt[4]{1-t}\right) - i\log\left(1-i\sqrt[4]{1-t}\right)\log\left(\sqrt[4]{1-t}\right) - iLi_2\left(i\sqrt[4]{1-t}\right)$$

(how could a human being deduce that?). We then have the definite integral

$$\frac{1}{4}\int_0^1 \frac{\log(1-t)}{t(1-t)^{3/4}} dt = Li_2(-1) - Li_2(1) + i\left[Li_2(i) - Li_2(-i)\right]$$

and reference to the definition of the dilogarithm function readily shows us that

$$Li_2(i) - Li_2(-i) = 2iG$$

and we therefore deduce that

(38.3) $$\varsigma\left(2,\frac{1}{4}\right) = \psi'\left(\frac{1}{4}\right) = \pi^2 + 8G$$

as previously derived by Kölbig [30i] in a much simpler manner.

In passing, it may be noted that whilst the beta function $B(u,v)$ [4] is only defined for $\text{Re}(u) > 0$ and $\text{Re}(v) > 0$, we may still determine $\left.\frac{\partial}{\partial v} B(u,v)\right|_{u=0}$ as follows:

Since $B(u,v) = \int_0^1 t^{u-1}(1-t)^{v-1} dt$ we have

$$\frac{\partial}{\partial v} B(u,v) = \int_0^1 t^{u-1}(1-t)^{v-1} \log(1-t) dt$$

and $B(u,v) = \frac{\Gamma(u)\Gamma(v)}{\Gamma(u+v)}$ also implies that

$$\frac{\partial}{\partial v} B(u,v) = B(u,v)\left[\psi(v) - \psi(u+v)\right]$$

$$= uB(u,v)\frac{\left[\psi(u+v) - \psi(v)\right]}{u}$$

and since $uB(u,v) = \dfrac{\Gamma(u+1)\Gamma(v)}{\Gamma(u+v)}$ we see that

$$\lim_{u \to 0} \frac{\partial}{\partial v} B(u,v) = -\lim_{u \to 0} \frac{\Gamma(u+1)\Gamma(v)}{\Gamma(u+v)} \lim_{u \to 0} \frac{[\psi(u+v) - \psi(v)]}{u}$$

$$= -\psi'(v)$$

Therefore we obtain

(38.4) $$\psi'(v) = -\int_0^1 \frac{(1-t)^{v-1} \log(1-t)}{t} dt$$

in agreement with (38.2). Equation (38.4) may also be derived by differentiating the well-known integral [56, p.15] for the digamma function

$$\psi(v) = -\gamma + \int_0^1 \frac{1 - t^{v-1}}{1-t} dt$$

From the definition of the Hurwitz zeta function we see that

$$\varsigma\left(2, \frac{1}{4}\right) - \varsigma\left(2, \frac{3}{4}\right) = 16 \sum_{n=0}^{\infty} \frac{1}{(4n+1)^2} - 16 \sum_{n=0}^{\infty} \frac{1}{(4n+3)^2}$$

$$= 16 \sum_{n=0}^{\infty} \left[ \frac{1}{(4n+1)^2} - \frac{1}{(4n+3)^2} \right] = 16G$$

For convenience we now designate $g(x)$ as

(38.5) $$g(x) = \frac{\Gamma(n+1)\Gamma(x)}{\Gamma(n+1+x)}$$

and differentiation gives us

$$g'(x) = g(x)[\psi(x) - \psi(n+1+x)]$$

Therefore differentiation of (34) results in

$$\varsigma'(2,x) = -2\varsigma(3,x) = \sum_{n=0}^{\infty} \frac{1}{n+1} \frac{\Gamma(n+1)\Gamma(x)}{\Gamma(n+1+x)} [\psi(x) - \psi(n+1+x)]$$

and since [4, p.13]

$$(38.6) \quad \psi(n+1+x)-\psi(x)=\frac{1}{x}+\frac{1}{x+1}+\ldots+\frac{1}{x+n}=H_{n+1}^{(1)}(x)$$

we obtain

$$(39) \quad 2\varsigma(3,x)=\sum_{n=0}^{\infty}\frac{1}{n+1}\frac{\Gamma(n+1)\Gamma(x)}{\Gamma(n+1+x)}H_{n+1}^{(1)}(x)=\sum_{n=1}^{\infty}\frac{1}{n}\frac{\Gamma(n)\Gamma(x)}{\Gamma(n+x)}H_n^{(1)}(x)$$

$$=\sum_{n=1}^{\infty}\frac{B(n,x)}{n}H_n^{(1)}(x)$$

With $x=1$ in (39) this becomes the very familiar sum originally derived by Euler

$$2\varsigma(3)=\sum_{n=0}^{\infty}\frac{H_{n+1}^{(1)}}{(n+1)^2}=\sum_{n=1}^{\infty}\frac{H_n^{(1)}}{n^2}$$

We note that

$$g'(x)=-g(x)H_{n+1}^{(1)}(x)$$

Differentiation of (39) results in

$$(40) \quad 3!\varsigma(4,x)=\sum_{n=1}^{\infty}\frac{1}{n}\frac{\Gamma(n)\Gamma(x)}{\Gamma(n+x)}\left(\left[H_n^{(1)}(x)\right]^2+H_n^{(2)}(x)\right)$$

and with $x=1$ we obtain a well-known result

$$(41) \quad \varsigma(4)=\frac{1}{3!}\sum_{n=1}^{\infty}\frac{\left[H_n^{(1)}\right]^2+H_n^{(2)}}{n^2}$$

which was also obtained in a different way in equation (4.2.40) of [12].

Differentiation of (40) gives us a more complex summation

$$(42) \quad \varsigma(5,x)=\frac{1}{4!}\sum_{n=1}^{\infty}\frac{1}{n}\frac{\Gamma(n)\Gamma(x)}{\Gamma(n+x)}\left(\left[H_n^{(1)}(x)\right]^3+3H_n^{(2)}(x)H_n^{(1)}(x)+2H_n^{(3)}(x)\right)$$

With $x=1$ we get another known result

(43) $$\varsigma(5) = \frac{1}{4!}\sum_{n=1}^{\infty}\frac{\left[H_n^{(1)}\right]^3 + 3H_n^{(2)}H_n^{(1)} + 2H_n^{(3)}}{n^2}$$

which was also previously derived in equation (4.2.47a) of [12] (and see the detailed references contained therein).

A further differentiation gives us

(43.1)

$$\varsigma(6,x) = \frac{1}{5!}\sum_{n=1}^{\infty}\frac{1}{n}\frac{\Gamma(n)\Gamma(x)}{\Gamma(n+x)}\left(3\left[H_n^{(1)}(x)\right]^2 H_n^{(2)}(x) + 3\left[H_n^{(2)}(x)\right]^2 + 3H_n^{(1)}(x)H_n^{(3)}(x) + 2H_n^{(4)}(x)\right)$$

$$+ \frac{1}{5!}\sum_{n=1}^{\infty}\frac{1}{n}\frac{\Gamma(n)\Gamma(x)}{\Gamma(n+x)}\left(\left[H_n^{(1)}(x)\right]^4 + 3\left[H_n^{(1)}(x)\right]^2 H_n^{(2)}(x) + 2H_n^{(1)}(x)H_n^{(3)}(x)\right)$$

and hence we have a further series involving non-linear Euler sums

(43.2) $$\varsigma(6) = \frac{1}{5!}\sum_{n=1}^{\infty}\frac{1}{n^2}\left(6\left[H_n^{(1)}\right]^2 H_n^{(2)} + 3\left[H_n^{(2)}\right]^2 + 5H_n^{(1)}H_n^{(3)} + 2H_n^{(4)} + \left[H_n^{(1)}\right]^4\right)$$

which we may compare with (45.10) below,

This differentiation procedure could obviously be extended to infinity and beyond! Instead, we now consider the Hasse/Coppo formula for $q = 2$. This gives us

$$\varsigma(3,x) = \frac{1}{2!}\sum_{n=1}^{\infty}\frac{1}{n}\frac{\Gamma(n)\Gamma(x)}{\Gamma(n+x)}H_n^{(1)}(x)$$

and this simply repeats (39) above (as indeed will also happen for $q = 3$ etc. due to the fact that (33.2) is the derivative of (33.1) and so on).

□

We note from (29.2) that

$$S(n,0,x-t) = S(n,0,x)\exp\left[\sum_{m=1}^{\infty}(m-1)!H_n^{(m)}(x)\frac{t^m}{m!}\right]$$

and by definition we have

$$S(n,0,x) = \frac{n!}{x(x+1)...(x+n-1)}$$

We have

$$\frac{n!}{x(x+1)\ldots(x+n-1)(x+n)} = \frac{\Gamma(n+1)\Gamma(x)}{\Gamma(n+1+x)}$$

and therefore we see that

$$S(n,0,x) = \frac{(x+n)\Gamma(n+1)\Gamma(x)}{\Gamma(n+1+x)} = \frac{\Gamma(n+1)\Gamma(x)}{\Gamma(n+x)}$$

We may therefore write (29.2) as

$$\frac{\Gamma(n+1)\Gamma(x-t)}{\Gamma(n+x-t)} = \frac{\Gamma(n+1)\Gamma(x)}{\Gamma(n+x)} \exp\left[\sum_{m=1}^{\infty} H_n^{(m)}(x) \frac{t^m}{m}\right]$$

or equivalently

$$\log \frac{\Gamma(n+x)\Gamma(x-t)}{\Gamma(x)\Gamma(n+x-t)} = \sum_{m=1}^{\infty} H_n^{(m)}(x) \frac{t^m}{m}$$

Let us now consider the case where $t \to -t$ and $x = 1$. This gives us

$$\log \frac{\Gamma(n+1)\Gamma(1+t)}{\Gamma(n+1+t)} = \sum_{m=1}^{\infty} \frac{(-1)^m}{m} H_n^{(m)} t^m$$

It is easily seen that

$$\log \frac{\Gamma(n+1)\Gamma(1+t)}{\Gamma(n+1+t)} = \log \frac{n\Gamma(n)\Gamma(1+t)}{(n+t)\Gamma(n+t)}$$

$$= \log \frac{\Gamma(n)\Gamma(1+t)}{\Gamma(n+t)} - \log\left(1 + \frac{t}{n}\right)$$

and we therefore have

$$\log \frac{\Gamma(n)\Gamma(1+t)}{\Gamma(n+t)} = \log\left(1 + \frac{t}{n}\right) + \sum_{m=1}^{\infty} \frac{(-1)^m}{m} H_n^{(m)} t^m$$

$$= -\sum_{m=1}^{\infty} \frac{(-1)^m}{m} \frac{t^m}{n^m} + \sum_{m=1}^{\infty} \frac{(-1)^m}{m} H_n^{(m)} t^m$$

$$= \sum_{m=1}^{\infty} \frac{(-1)^m}{m} \left( H_n^{(m)} - \frac{1}{n^m} \right) t^m$$

We then end up with

(44) $$\log \frac{\Gamma(n)\Gamma(1+t)}{\Gamma(n+t)} = \sum_{m=1}^{\infty} \frac{(-1)^m}{m} H_{n-1}^{(m)} t^m$$

## DETERMINATION OF THE STIRLING NUMBERS $s(n,k)$

Some time ago the author considered a similar expression to (44) in [13]. Let

$$f(x) = \log \Gamma(n+x) - \log \Gamma(1+x)$$

Then we have

$$f^{(m)}(0) = \psi^{(m-1)}(n) - \psi^{(m-1)}(1) = (-1)^{m+1}(m-1)! H_{n-1}^{(m)}$$

and we therefore have the Maclaurin expansion equivalent to (44)

(44.1) $$f(x) = \log \Gamma(n+x) - \log \Gamma(1+x) = \log \Gamma(n) + \sum_{m=1}^{\infty} \frac{(-1)^{m+1}}{m} H_{n-1}^{(m)} x^m$$

This may be written as

(44.2) $$\frac{\Gamma(n+x)}{\Gamma(1+x)\Gamma(n)} = \exp\left[\sum_{m=1}^{\infty} \frac{(-1)^{m+1}}{m} H_{n-1}^{(m)} x^m\right]$$

and expansion of the exponential function leads to

(44.3)

$$\frac{\Gamma(n+x)}{\Gamma(1+x)\Gamma(n)} = 1 + H_{n-1}^{(1)} x + \frac{1}{2}\left(\left[H_{n-1}^{(1)}\right]^2 - H_{n-1}^{(2)}\right) x^2 + \frac{1}{6}\left(\left[H_{n-1}^{(1)}\right]^3 - 3 H_{n-1}^{(1)} H_{n-1}^{(2)} + 2 H_{n-1}^{(3)}\right) x^3 + O(x^4)$$

and

(44.4)

$$\frac{\Gamma(1+x)\Gamma(n)}{\Gamma(n+x)} = 1 - H_{n-1}^{(1)} x + \frac{1}{2}\left(\left[H_{n-1}^{(1)}\right]^2 + H_{n-1}^{(2)}\right) x^2 - \frac{1}{6}\left(\left[H_{n-1}^{(1)}\right]^3 + 3 H_{n-1}^{(1)} H_{n-1}^{(2)} + 2 H_{n-1}^{(3)}\right) x^3 + O(x^4)$$

The terms involving the generalised harmonic numbers are truly ubiquitous! I came across (44.2) in a 1999 paper entitled "Analytic two-loop results for self energy- and vertex-type diagrams with one non-zero mass" by Fleisher et al [20]. Letting $n = 2$ in (44.1) results in the familiar Maclaurin expansion for $\log(1+x)$.

From (44.1) and using the formulae (20.1) for the Stirling numbers of the first kind we see that

$$(-1)^{n+1} \frac{\Gamma(n+x)}{\Gamma(1+x)} = s(n,1) - s(n,2)x + s(n,3)x^2 - s(n,4)x^3 + O(x^4)$$

and this may be written using the Pochhammer symbol $(x)_n$ as

(44.5) $$(x)_n = \frac{\Gamma(n+x)}{\Gamma(x)} = \frac{x\Gamma(n+x)}{\Gamma(1+x)} = \sum_{k=0}^{\infty} (-1)^{n+k} s(n,k) x^k$$

which is simply the generating function for the Stirling numbers of the first kind (letting $x \to -x$ in (20.1))

$$(x)_n = x(x+1)...(x+n-1) = \sum_{k=0}^{\infty} (-1)^{n+k} s(n,k) x^k = \sum_{k=1}^{\infty} (-1)^{n+k} s(n,k) x^k$$

since $s(n,k) = 0$ for $k \geq n+1$. Using (44.5) is probably the simplest way of evaluating the Stirling numbers (by successive differentiation); for example we have

(44.6) $$\frac{\Gamma(n+x)}{\Gamma(1+x)} \left[ x\psi(n+x) + 1 - x\psi(1+x) \right] = \sum_{k=1}^{\infty} (-1)^{n+k} k s(n,k) x^{k-1}$$

and with $x = 0$ we obtain $s(n,1) = (-1)^{n+1}(n-1)!$.

A further differentiation results in

$$\frac{\Gamma(n+x)}{\Gamma(1+x)} \left[ \psi(n+x) - \psi(1+x) \right] \left[ x\psi(n+x) + 1 - x\psi(1+x) \right]$$

$$+ \frac{\Gamma(n+x)}{\Gamma(1+x)} \left( x \left[ \psi'(n+x) - \psi'(1+x) \right] + \left[ \psi(n+x) - \psi(1+x) \right] \right) = \sum_{k=2}^{\infty} (-1)^{n+k} k(k-1) s(n,k) x^{k-2}$$

and with $x = 0$ we easily obtain $s(n,2) = (-1)^n (n-1)! H_{n-1}^{(1)}$.

These results may be generalised as follows. The Pochhammer symbol may be expressed as

$$(u+x)_n = (u+x)(1+u+x)\ldots(n-1+u+x)$$

$$= \left(1+\frac{x}{u}\right)\left(1+\frac{x}{1+u}\right)\ldots\left(1+\frac{x}{n-1+u}\right)u(1+u)(2+u)\ldots(n-1+u)$$

$$= \left(1+\frac{x}{u}\right)\left(1+\frac{x}{1+u}\right)\ldots\left(1+\frac{x}{n-1+u}\right)(u)_n$$

and we see that

$$\left(1+\frac{x}{u}\right)\left(1+\frac{x}{1+u}\right)\ldots\left(1+\frac{x}{n-1+u}\right) = \exp\left[\sum_{j=0}^{n-1}\log\left(1+\frac{x}{j+u}\right)\right]$$

$$= \exp\left[\sum_{j=0}^{n-1}\sum_{m=1}^{\infty}(-1)^{m-1}\frac{x^m}{m(j+u)^m}\right]$$

$$= \exp\left[\sum_{m=1}^{\infty}(-1)^{m-1}\frac{x^m}{m}\sum_{j=0}^{n-1}\frac{1}{(j+u)^m}\right]$$

$$= \exp\left[\sum_{m=1}^{\infty}(-1)^{m-1}\frac{x^m}{m}H_n^{(m)}(u)\right]$$

We therefore obtain

$$\frac{(u+x)_n}{(u)_n} = \exp\left[\sum_{m=1}^{\infty}(-1)^{m-1}\frac{x^m}{m}H_n^{(m)}(u)\right]$$

and we know that

$$(u+x)_n = \sum_{k=0}^{n}(-1)^{n+k}s(n,k)(u+x)^k$$

$$= \sum_{k=0}^{n}(-1)^{n+k}s(n,k)\sum_{j=0}^{k}\binom{k}{j}u^{k-j}x^j$$

$$= \sum_{j=0}^{n}\sum_{k=j}^{n}(-1)^{n+k}s(n,k)\binom{k}{j}u^{k-j}x^j$$

Differentiation results in

$$\frac{d^r}{dx^r}(u+x)_n = \sum_{j=0}^{n}\sum_{k=j}^{n}(-1)^{n+k}s(n,k)\binom{k}{j}u^{k-j}j(j-1)..(j-r+1)x^{j-r}$$

and we see that

$$\left.\frac{d^r}{dx^r}(u+x)_n\right|_{x=0} = \sum_{k=r}^{n}(-1)^{n+k}s(n,k)r!\binom{k}{r}u^{k-r}$$

We also have from (24)

$$\frac{d^r}{dx^r}\exp\left[\sum_{m=1}^{\infty}(-1)^{m+1}\frac{x^m}{m}H_n^{(m)}(u)\right] = e^{f(x)}Y_r\left(f^{(1)}(x), f^{(2)}(x), \ldots, f^{(r)}(x)\right)$$

where $f(x) = \sum_{m=1}^{\infty}(-1)^{m+1}\frac{x^m}{m}H_n^{(m)}(u)$ and we have

$$f^{(p)}(x) = \sum_{m=1}^{\infty}(-1)^{m+1}m(m-1)\ldots(m-p+1)\frac{x^{m-p}}{m}H_n^{(m)}(u)$$

With $x = 0$ we see that

$$f^{(p)}(0) = (-1)^{p-1}(p-1)!H_n^{(p)}(u)$$

and we thereby obtain

$$\left.\frac{d^r}{dx^r}\exp\left[\sum_{m=1}^{\infty}(-1)^{m+1}\frac{x^m}{m}H_n^{(m)}(u)\right]\right|_{x=0} = (u)_n Y_r\left(H_n^{(1)}(u), -1!H_n^{(2)}(u), \ldots, (-1)^{r-1}(r-1)!H_n^{(r)}(u)\right)$$

This then gives us

$$(44.7)\quad r!\sum_{k=r}^{n}(-1)^{n+k}s(n,k)\binom{k}{r}u^{k-r} = (u)_n Y_r\left(H_n^{(1)}(u), -1!H_n^{(2)}(u), \ldots, (-1)^{r-1}(r-1)!H_n^{(r)}(u)\right)$$

and in particular, with $u = 1$, we obtain using $(1)_n = n!$

$$(44.8)\quad \sum_{k=r}^{n}(-1)^{n+k}s(n,k)\binom{k}{r} = \frac{n!}{r!}Y_r\left(H_n^{(1)}, -1!H_n^{(2)}, \ldots, (-1)^{r-1}(r-1)!H_n^{(r)}\right)$$

Since $(x)_{n+1} = x(1+x)_n$ we can show that

(44.9) $$(-1)^{n+r} s(n+1, r+1) = \sum_{k=r}^{n} (-1)^{n+k} s(n,k) \binom{k}{r}$$

and hence, as reported by Comtet [10], we have in terms of the complete Bell polynomials

(44.10) $$s(n+1, r+1) = (-1)^{n+r} \frac{n!}{r!} Y_r\left(H_n^{(1)}, -1! H_n^{(2)}, \ldots, (-1)^{r-1}(r-1)! H_n^{(r)}\right)$$

I subsequently discovered that this is a slightly different version of the proof originally given by Kölbig [30] in his 1993 thesis.

## FURTHER EXAMPLES OF EULER-HURWITZ SUMS

Employing the Nörlund representation for the digamma function (see equation (36) above), it was also shown in equation (4.3.48) in [12] that

(45) $$\psi^{(q)}(x) = q! \sum_{n=1}^{\infty} \frac{(-1)^{n+1}}{n} \frac{s(n,q)}{x(x+1)\ldots(x+n-1)}$$

and with $x = 1$ we obtain

(45.1) $$\psi^{(q)}(1) = q! \sum_{n=1}^{\infty} \frac{(-1)^{n+1}}{n \cdot n!} s(n,q)$$

Since [46, p.22] $\psi^{(q)}(x) = (-1)^{q+1} q! \varsigma(q+1, x)$ we have

(45.2) $$\varsigma(p+1) = (-1)^p \sum_{k=1}^{\infty} \frac{(-1)^k}{k \cdot k!} s(k, p)$$

This result was previously obtained in 1995 by Shen [41] by employing a different method. Other proofs were recently given in [7] and [46, p.252].

We may therefore write (45) as

$$\varsigma(q+1, x) = (-1)^{q+1} \sum_{n=1}^{\infty} \frac{(-1)^{n+1}}{n} \frac{s(n,q)}{x(x+1)\ldots(x+n-1)}$$

$$= (-1)^{q+1} \sum_{n=1}^{\infty} \frac{(-1)^{n+1}}{n} \frac{s(n,q)(n+x)\Gamma(x)}{\Gamma(n+1+x)}$$

and therefore we obtain

$$(45.3) \qquad \varsigma(q+1,x) = (-1)^{q+1} \sum_{n=1}^{\infty} \frac{(-1)^{n+1}}{n} \frac{s(n,q)\Gamma(x)}{\Gamma(n+x)}$$

We may write (44.10) as

$$s(n,q) = (-1)^{n+q} \frac{(n-1)!}{(q-1)!} Y_{q-1}\left(H_{n-1}^{(1)}, -1! H_{n-1}^{(2)}, ..., (-1)^q (q-2)! H_{n-1}^{(q-1)}\right)$$

and substituting this in (45.3) gives us

$$\varsigma(q+1,x) = (-1)^{q+1} \sum_{n=1}^{\infty} \frac{(-1)^{n+1}}{n} \frac{\Gamma(x)}{\Gamma(n+x)} (-1)^{n+q} \frac{(n-1)!}{(q-1)!} Y_{q-1}\left(H_{n-1}^{(1)}, -1! H_{n-1}^{(2)}, ..., (-1)^q (q-2)! H_{n-1}^{(q-1)}\right)$$

This may be slightly simplified to

$$(45.4) \qquad \varsigma(q+1,x) = \frac{1}{(q-1)!} \sum_{n=1}^{\infty} \frac{1}{n} \frac{\Gamma(n)\Gamma(x)}{\Gamma(n+x)} Y_{q-1}\left(H_{n-1}^{(1)}, -1! H_{n-1}^{(2)}, ..., (-1)^q (q-2)! H_{n-1}^{(q-1)}\right)$$

where we particularly note that in this expression the complete Bell polynomials do **not** contain the $x$ variable (and is therefore a completely different representation from (32) which is reproduced below for ease of comparison).

$$\varsigma(q+1,x) = \frac{1}{q!} \sum_{n=1}^{\infty} \frac{1}{n} \frac{\Gamma(n)\Gamma(x)}{\Gamma(n+x)} Y_{q-1}\left(0! H_n^{(1)}(x), 1! H_n^{(2)}(x), ..., (q-2)! H_n^{(q-1)}(x)\right)$$

With $q=1$ in (45.4) we get

$$\varsigma(2,x) = \sum_{n=1}^{\infty} \frac{1}{n} \frac{\Gamma(n)\Gamma(x)}{\Gamma(n+x)} Y_0 = \sum_{n=1}^{\infty} \frac{1}{n} \frac{\Gamma(n)\Gamma(x)}{\Gamma(n+x)}$$

and we have seen this before in (34).

Letting $q=2$ gives us

$$\varsigma(3,x) = \sum_{n=1}^{\infty} \frac{1}{n} \frac{\Gamma(n)\Gamma(x)}{\Gamma(n+x)} Y_1\left(H_{n-1}^{(1)}\right) = \sum_{n=1}^{\infty} \frac{1}{n} \frac{\Gamma(n)\Gamma(x)}{\Gamma(n+x)} H_{n-1}^{(1)}$$

$$= \sum_{n=1}^{\infty} \frac{1}{n} \frac{\Gamma(n)\Gamma(x)}{\Gamma(n+x)} \left(H_n^{(1)} - \frac{1}{n}\right)$$

We then have

$$(45.5) \qquad \varsigma(3,x) = \sum_{n=1}^{\infty} \frac{1}{n} \frac{\Gamma(n)\Gamma(x)}{\Gamma(n+x)} H_n^{(1)} - \sum_{n=1}^{\infty} \frac{1}{n^2} \frac{\Gamma(n)\Gamma(x)}{\Gamma(n+x)}$$

and this may be contrasted with (39)

$$\varsigma(3,x) = \frac{1}{2} \sum_{n=1}^{\infty} \frac{1}{n} \frac{\Gamma(n)\Gamma(x)}{\Gamma(n+x)} H_n^{(1)}(x)$$

Letting $x=1$ in (45.5) reproduces the Euler sum

$$\varsigma(3) = \sum_{n=1}^{\infty} \frac{H_n^{(1)}}{n^2} - \sum_{n=1}^{\infty} \frac{1}{n^3}$$

and with $x = 1/2$ we obtain using (36.2)

$$(45.6) \quad \varsigma\left(3,\frac{1}{2}\right) = \sum_{n=1}^{\infty} \frac{H_n^{(1)}\left(\frac{1}{2}\right)}{n} \frac{2^{2n-1}\Gamma^2(n)}{\Gamma(2n)} - \sum_{n=1}^{\infty} \frac{1}{n^2} \frac{2^{2n-1}\Gamma^2(n)}{\Gamma(2n)} = \frac{1}{2}\sum_{n=1}^{\infty} \frac{nH_n^{(1)}\left(\frac{1}{2}\right)-1}{n^2} \frac{\left[2^n\Gamma(n)\right]^2}{\Gamma(2n)}$$

where $H_n^{(1)}\left(\frac{1}{2}\right) = 2 \sum_{k=0}^{n-1} \frac{1}{2k+1}$.

Differentiation of (45.5) results in

$$\varsigma(4,x) = \frac{1}{3}\sum_{n=1}^{\infty} \frac{H_n^{(1)}}{n} \frac{\Gamma(n)\Gamma(x)}{\Gamma(n+x)}[\psi(n+x)-\psi(x)] - \frac{1}{3}\sum_{n=1}^{\infty} \frac{1}{n^2} \frac{\Gamma(n)\Gamma(x)}{\Gamma(n+x)}[\psi(n+x)-\psi(x)]$$

Using (38.6) we see that

$$\psi(n+x) - \psi(x) = H_n^{(1)}(x)$$

We therefore have

$$\varsigma(4,x) = \frac{1}{3}\sum_{n=1}^{\infty} \frac{H_n^{(1)} H_n^{(1)}(x)}{n} \frac{\Gamma(n)\Gamma(x)}{\Gamma(n+x)} - \frac{1}{3}\sum_{n=1}^{\infty} \frac{H_n^{(1)}(x)}{n^2} \frac{\Gamma(n)\Gamma(x)}{\Gamma(n+x)}$$

or

$$(45.7) \qquad \varsigma(4,x) = \frac{1}{3}\sum_{n=1}^{\infty} \frac{\left[nH_n^{(1)}-1\right]H_n^{(1)}(x)}{n^2} \frac{\Gamma(n)\Gamma(x)}{\Gamma(n+x)}$$

Letting $x = 1$ we get

$$\varsigma(4) = \frac{1}{3}\sum_{n=1}^{\infty} \frac{\left[H_n^{(1)}\right]^2}{n^2} - \frac{1}{3}\sum_{n=1}^{\infty} \frac{H_n^{(1)}}{n^3}$$

and this concurs with the following well-known results (see, for example, equation (4.3.46b) in [12])

$$\sum_{n=1}^{\infty} \frac{H_n^{(1)}}{n^3} = \frac{5}{4}\varsigma(4) \qquad \sum_{n=1}^{\infty} \frac{\left[H_n^{(1)}\right]^2}{n^2} = \frac{17}{4}\varsigma(4)$$

It is known from [34, p.279] that

$$nH_n^{(1)} = n + \sum_{k=1}^{n-1} H_k^{(1)}$$

and therefore (45.7) may be written as

$$\varsigma(4, x) = \frac{1}{3}\sum_{n=1}^{\infty} \frac{\left[n - 1 + \sum_{k=1}^{n-1} H_k^{(1)}\right] H_n^{(1)}(x)}{n^2} \frac{\Gamma(n)\Gamma(x)}{\Gamma(n+x)}$$

Differentiation of (45.7) results in

$$\varsigma(5, x) = \frac{1}{12}\sum_{n=1}^{\infty} \frac{\left[nH_n^{(1)} - 1\right]\left[H_n^{(1)}(x)\right]^2}{n^2} \frac{\Gamma(n)\Gamma(x)}{\Gamma(n+x)} + \frac{1}{12}\sum_{n=1}^{\infty} \frac{\left[nH_n^{(1)} - 1\right] H_n^{(2)}(x)}{n^2} \frac{\Gamma(n)\Gamma(x)}{\Gamma(n+x)}$$

which may be written as

(45.7) $$\varsigma(5, x) = \frac{2}{4!}\sum_{n=1}^{\infty} \frac{\left[nH_n^{(1)} - 1\right]\left(\left[H_n^{(1)}(x)\right]^2 + H_n^{(2)}(x)\right)}{n^2} \frac{\Gamma(n)\Gamma(x)}{\Gamma(n+x)}$$

We then have with $x = 1$

$$\varsigma(5) = \frac{1}{12}\sum_{n=1}^{\infty} \frac{\left[nH_n^{(1)} - 1\right]\left(\left[H_n^{(1)}\right]^2 + H_n^{(2)}\right)}{n^3}$$

$$= \frac{1}{12}\sum_{n=1}^{\infty}\frac{\left[H_n^{(1)}\right]^3 + H_n^{(1)}H_n^{(2)}}{n^2} - \frac{1}{12}\sum_{n=1}^{\infty}\frac{\left[H_n^{(1)}\right]^2 + H_n^{(2)}}{n^3}$$

and therefore

(45.8) $$12\varsigma(5) = \sum_{n=1}^{\infty}\frac{\left[H_n^{(1)}\right]^3}{n^2} + \sum_{n=1}^{\infty}\frac{H_n^{(1)}H_n^{(2)}}{n^2} - \sum_{n=1}^{\infty}\frac{\left[H_n^{(1)}\right]^2}{n^3} - \sum_{n=1}^{\infty}\frac{H_n^{(2)}}{n^3}$$

and this concurs with the four well-known individual Euler sums (which are also derived in [12]).

Differentiating (45.7) gives us

(45.9) $$\varsigma(6,x) = \frac{1}{60}\sum_{n=1}^{\infty}\frac{\left[nH_n^{(1)} - 1\right]\left(\left[H_n^{(1)}(x)\right]^2 + H_n^{(2)}(x)\right)}{n^2}\frac{\Gamma(n)\Gamma(x)}{\Gamma(n+x)}H_n^{(1)}(x)$$

$$+ \frac{1}{60}\sum_{n=1}^{\infty}\frac{\left[nH_n^{(1)} - 1\right]\left(2H_n^{(1)}(x)H_n^{(2)}(x) + H_n^{(3)}(x)\right)}{n^2}\frac{\Gamma(n)\Gamma(x)}{\Gamma(n+x)}$$

and with $x = 1$ we obtain

(45.10) $$\frac{1}{2}5!\varsigma(6) = \sum_{n=1}^{\infty}\frac{\left[H_n^{(1)}\right]^4 + 3\left[H_n^{(1)}\right]^2 H_n^{(2)} + H_n^{(1)}H_n^{(3)}}{n^2} - \sum_{n=1}^{\infty}\frac{\left[H_n^{(1)}\right]^3 + 3H_n^{(1)}H_n^{(2)} + H_n^{(3)}}{n^3}$$

which may be written as

$$\frac{1}{2}5!\varsigma(6) = \sum_{n=1}^{\infty}\frac{\left(nH_n^{(1)} - 1\right)\left(\left[H_n^{(1)}\right]^3 + 3H_n^{(1)}H_n^{(2)} + H_n^{(3)}\right)}{n^3}$$

Four other representations of $\varsigma(6)$ involving the generalised harmonic numbers are derived in equation (4.3.65) in [12]. It is possible that this set of simultaneous linear equations, in conjunction with Zheng's recent paper [51] and other known results, may be sufficient to evaluate more complex non-linear Euler sums and perhaps verify Coffey's conjecture [7i] that

$$\sum_{n=1}^{\infty}\frac{\left[H_n^{(1)}\right]^4}{(n+1)^2} = \frac{859}{24}\varsigma(6) + 3\varsigma^2(3)$$

which we may express in standard form by writing $H_n^{(1)} = H_{n+1}^{(1)} - \frac{1}{n+1}$.

## A CONNECTION WITH THE GAUSS REPRESENTATION OF THE GAMMA FUNCTION

We recall (44.2) and, as noted by Wilf [50], we have

$$\frac{\Gamma(n+x)}{\Gamma(1+x)\Gamma(n)} = \exp\left[\sum_{m=1}^{\infty} \frac{(-1)^{m+1}}{m} H_{n-1}^{(m)} x^m\right]$$

$$= \exp\left[H_{n-1}^{(1)} x\right] \exp\left[\sum_{m=2}^{\infty} \frac{(-1)^{m+1}}{m} H_{n-1}^{(m)} x^m\right]$$

$$= \exp\left[\left(H_{n-1}^{(1)} - \log n\right) x\right] \exp(x \log n) \exp\left[\sum_{m=2}^{\infty} \frac{(-1)^{m+1}}{m} H_{n-1}^{(m)} x^m\right]$$

We therefore have

$$\exp\left[\sum_{m=2}^{\infty} \frac{(-1)^{m+1}}{m} H_{n-1}^{(m)} x^m\right] = \frac{\Gamma(n+x)}{\Gamma(1+x)\Gamma(n) n^x \exp\left[\left(H_{n-1}^{(1)} - \log n\right) x\right]}$$

and, using (17.2), this may be expressed as

$$\exp\left[\sum_{m=2}^{\infty} \frac{(-1)^{m+1}}{m} H_{n-1}^{(m)} x^m\right] = \frac{1}{x} \frac{n}{n+x} \frac{x(1+x)\ldots(n+x)}{n! n^x} \frac{1}{\exp\left[\left(H_{n-1}^{(1)} - \log n\right) x\right]}$$

We note from (9) that

$$\lim_{n \to \infty} \exp\left[\left(H_{n-1}^{(1)} - \log n\right) x\right] = e^{\gamma x}$$

and we see that

$$\lim_{n \to \infty} \exp\left[\sum_{m=2}^{\infty} \frac{(-1)^{m+1}}{m} H_{n-1}^{(m)} x^m\right] = \exp\left[\sum_{m=2}^{\infty} \frac{(-1)^{m+1}}{m} \varsigma(m) x^m\right]$$

We have the Gauss expression [46, p.2] for the gamma function

$$\Gamma(x) = \lim_{n \to \infty}\left[\frac{n!n^x}{x(1+x)\ldots(n+x)}\right]$$

and hence we obtain

$$\exp\left[\sum_{m=2}^{\infty}\frac{(-1)^{m+1}}{m}\varsigma(m)x^m\right] = \frac{e^{-\gamma x}}{\Gamma(1+x)}$$

which is equivalent to the well known series expansion

(46) $$\log\Gamma(1+x) = -\gamma x + \sum_{m=2}^{\infty}\frac{(-1)^{m+1}}{m}\varsigma(m)x^m$$

With regard to (24) we let

$$f(x) = \log(u+x)_n$$

$$= \sum_{j=0}^{n-1}\log(j+u+x)$$

and therefore we see that

$$f^{(p)}(x) = (-1)^{p-1}(p-1)!\sum_{j=0}^{n-1}\frac{1}{(j+u+x)^p} = (-1)^{p-1}(p-1)!H_n^{(p)}(u+x)$$

Hence from (24) we get as before

$$\frac{d^m}{dx^m}(u+x)_n = (u+x)_n Y_m\left(H_n^{(1)}(u+x), -1!H_n^{(2)}(u+x), \ldots, (-1)^{m-1}(m-1)!H_n^{(m)}(u+x)\right)$$

or equivalently by reference to (44.5)

$$\frac{d^m}{dx^m}\frac{\Gamma(n+u+x)}{\Gamma(u+x)} = \frac{\Gamma(n+u+x)}{\Gamma(u+x)}Y_m\left(H_n^{(1)}(u+x), -1!H_n^{(2)}(u+x), \ldots, (-1)^{m-1}(m-1)!H_n^{(m)}(u+x)\right)$$

## SOME INTEGRALS

We have

$$\frac{\Gamma(1+x)\Gamma(n)}{\Gamma(n+x)} = \frac{x\Gamma(x)\Gamma(n)}{\Gamma(n+x)} = xB(x,n)$$

It is an exercise in Whittaker & Watson [49, p.262] to show that

(47) $$\log B(x,n) = \log\left(\frac{x+n}{xn}\right) + \int_0^1 \frac{(1-v^x)(1-v^n)}{(1-v)\log v} dv \qquad x,n > 0$$

and this formula is attributed to Euler (see also [34, p.187]). In passing, we note the required symmetry in $x$ and $n$.

Therefore we have using (44.2)

(47.1) $$-\sum_{k=1}^{\infty} \frac{(-1)^{k+1}}{k} H_{n-1}^{(k)} x^k = \log\left(1+\frac{x}{n}\right) + \int_0^1 \frac{(1-v^x)(1-v^n)}{(1-v)\log v} dv$$

and

$$\log\frac{\Gamma(n+x)}{\Gamma(1+x)\Gamma(n)} = \log\left(1+\frac{x}{n}\right) + \int_0^1 \frac{(1-v^x)(1-v^n)}{(1-v)\log v} dv$$

An alternative proof is shown below.

We see that $1-v^x = 1-e^{x\log v} = -\sum_{k=1}^{\infty} \frac{x^k \log^k v}{k!}$ and therefore we have

$$\int_0^1 \frac{(1-v^x)(1-v^n)}{(1-v)\log v} dv = -\sum_{k=1}^{\infty} \frac{x^k}{k!} \int_0^1 \frac{(1-v^n)\log^{k-1} v}{1-v} dv$$

$$= -\sum_{k=1}^{\infty} \frac{x^k}{k!} \sum_{j=1}^{n-1} \int_0^1 v^j \log^{k-1} v \, dv$$

We have $\int_0^1 v^j v^a \, dv = \frac{1}{j+a+1}$ and therefore differentiation with respect to $a$ results in

$$\int_0^1 v^j v^a \log^{k-1} v \, dv = \frac{(-1)^{k-1}(k-1)!}{(j+a+1)^k}$$

Hence we have

$$\int_0^1 v^j \log^{k-1} v \, dv = \frac{(-1)^{k-1}(k-1)!}{(j+1)^k}$$

$$\sum_{j=0}^{n-1} \int_0^1 v^j \log^{k-1} v \, dv = (-1)^{k-1}(k-1)! H_n^{(k)}$$

and thus

(47.2) $$\int_0^1 \frac{(1-v^x)(1-v^n)}{(1-v)\log v}dv = -\sum_{k=1}^{\infty}\frac{(-1)^{k+1}}{k}H_n^{(k)}x^k$$

$$= -\sum_{k=1}^{\infty}\frac{(-1)^{k+1}}{k}\left(H_{n-1}^{(k)}+\frac{1}{n^k}\right)x^k$$

$$= -\sum_{k=1}^{\infty}\frac{(-1)^{k+1}}{k}H_{n-1}^{(k)}x^k - \sum_{k=1}^{\infty}\frac{(-1)^{k+1}}{k}\frac{x^k}{n^k}$$

$$= -\sum_{k=1}^{\infty}\frac{(-1)^{k+1}}{k}H_{n-1}^{(k)}x^k - \log\left(1+\frac{x}{n}\right)$$

Differentiating (47.2) results in

$$\int_0^1 \frac{v^x(1-v^n)}{1-v}dv = \sum_{k=1}^{\infty}(-1)^{k+1}H_n^{(k)}x^{k-1}$$

and as $x \to 0$ we have the well-known in integral

$$\int_0^1 \frac{(1-v^n)}{1-v}dv = H_n^{(1)}$$

Successive differentiations give us

$$\int_0^1 \frac{v^x(1-v^n)\log^p v}{1-v}dv = \sum_{k=p+1}^{\infty}(-1)^{k+1}H_n^{(k)}(k-1)(k-2)...(k-p)x^{k-p-1}$$

and hence with $x=1$ we see that

(47.3) $$\int_0^1 \frac{(1-v^n)\log^p v}{1-v}dv = (-1)^p p! H_n^{(p+1)}$$

We also have from equation (E.35c) in [17]

(47.4) $$\int_0^1 \frac{(1-v^n)v^{x-1}\log^p v}{1-v}dv = (-1)^p p!\sum_{k=1}^{n}\frac{1}{(k+x-1)^{p+1}} = (-1)^p p! H_n^{(p+1)}(x)$$

and this may also be derived by differentiating (38.4)

$$\int_0^1 \frac{(1-t)^{x-1} \log(1-t)}{t} dt = -\psi'(x)$$

Completing the summation of (47.4) gives us for $|u|<1$

$$(-1)^p p! \sum_{n=1}^\infty \frac{H_n^{(p+1)}(x)}{n} u^n = \sum_{n=1}^\infty \frac{1}{n} u^n \int_0^1 \frac{(1-v^n)v^{x-1}\log^p v}{1-v} dv$$

We therefore obtain

(47.5) $$(-1)^p p! \sum_{n=1}^\infty \frac{H_n^{(p+1)}(x)}{n} u^n = \int_0^1 \frac{[\log(1-uv) - \log(1-u)]v^{x-1}\log^p v}{1-v} dv$$

With $p=1$ we have

(47.6) $$\sum_{n=1}^\infty \frac{H_n^{(2)}(x)}{n} u^n = -\int_0^1 \frac{[\log(1-uv) - \log(1-u)]v^{x-1}\log v}{1-v} dv$$

and with $x=1$ we have

(47.7) $$\sum_{n=1}^\infty \frac{H_n^{(2)}}{n} u^n = -\int_0^1 \frac{[\log(1-uv) - \log(1-u)]\log v}{1-v} dv$$

The Wolfram Integrator is only able to evaluate (47.7) in terms involving the dilogarithm and the trilogarithm in the cases where $p = 0, 1$.

We note from equation (3.106f) in [11] that in terms of the polylogarithm function (14.5)

$$\sum_{n=1}^\infty \frac{H_n^{(2)}}{n} u^n = -2Li_3\left(-\frac{u}{1-u}\right)$$

Integrating (47.5) will give us an expression for $\sum_{n=1}^\infty \frac{H_n^{(p+1)}}{n^2} u^n$.

Using (47.4) we may make a more general summation for $\operatorname{Re}(s) > 1$

$$(-1)^p p! \sum_{n=1}^\infty \frac{H_n^{(p+1)}(x)}{n^s} u^n = \sum_{n=1}^\infty \frac{u^n}{n^s} \int_0^1 \frac{(1-v^n)v^{x-1}\log^p v}{1-v} dv$$

and therefore we get

(47.6) $$(-1)^p p! \sum_{n=1}^{\infty} \frac{H_n^{(p+1)}(x)}{n^s} u^n = \int_0^1 \frac{[Li_s(u) - Li_s(uv)] v^{x-1} \log^p v}{1-v} dv$$

We have a similar relationship from (4.4.156c) in [15]

$$\sum_{n=1}^{\infty} \frac{H_n^{(1)}}{n^s} u^n = \int_0^u \frac{Li_s(u) - Li_s(v)}{u - v} dv$$

With $u = 1$ in (47.6) we obtain

(47.7) $$(-1)^p p! \sum_{n=1}^{\infty} \frac{H_n^{(p+1)}(x)}{n^s} = \int_0^1 \frac{[\varsigma(s) - Li_s(v)] v^{x-1} \log^p v}{1-v} dv$$

and with $x = 1$

(47.8) $$(-1)^p p! \sum_{n=1}^{\infty} \frac{H_n^{(p+1)}}{n^s} = \int_0^1 \frac{[\varsigma(s) - Li_s(v)] \log^p v}{1-v} dv$$

□

In (39) we saw that

$$2\varsigma(3, x) = \sum_{n=1}^{\infty} \frac{B(n, x)}{n} H_n^{(1)}(x)$$

where $B(u, v)$ is Euler's beta function [4] defined for $\operatorname{Re}(u) > 0$ and $\operatorname{Re}(v) > 0$

$$B(u, v) = \int_0^1 t^{u-1} (1-t)^{v-1} dt \quad \text{and} \quad B(u, v) = \frac{\Gamma(u)\Gamma(v)}{\Gamma(u+v)}$$

From (47.4) we have

$$H_n^{(1)}(x) = \int_0^1 \frac{(1-v^n) v^{x-1}}{1-v} dv$$

and we therefore obtain a triple integral representation for $\varsigma(3, x)$

$$\varsigma(3, x) = \frac{1}{2} \sum_{n=1}^{\infty} \int_0^1 u^{n-1} du \int_0^1 t^{x-1} (1-t)^{n-1} dt \int_0^1 \frac{(1-v^n) v^{x-1}}{1-v} dv$$

The geometric series then gives us

$$\varsigma(3,x) = \frac{1}{2}\int_0^1\int_0^1\int_0^1 \frac{(tv)^{x-1}}{1-v}\left[\frac{1}{1-u(1-t)} - \frac{v}{1-uv(1-t)}\right] dt\, du\, dv$$

We have

$$\int_0^1 \frac{v}{1-uv(1-t)} du = -\frac{\log[1-v(1-t)]}{1-t}$$

and hence we obtain the double integral representation

$$\varsigma(3,x) = \frac{1}{2}\int_0^1\int_0^1 \frac{(tv)^{x-1}\left(\log[1-v(1-t)] - \log t\right)}{(1-t)(1-v)} dt\, dv$$

With $x=1$ we have

$$\varsigma(3) = \frac{1}{2}\int_0^1\int_0^1 \frac{\log[1-v(1-t)] - \log t}{(1-t)(1-v)} dt\, dv$$

AN INTEGRAL INVOLVING THE HURWITZ-LERCH ZETA FUNCTION

Differentiation of (47.1) gives us

(48) $$\psi(n+x) - \psi(1+x) = \frac{1}{n+x} - \int_0^1 \frac{v^x(1-v^n)}{1-v} dv$$

and with $x = 1/2$ we have

$$\psi\left(n+\frac{1}{2}\right) - \psi\left(\frac{3}{2}\right) = \frac{2}{2n+1} - \int_0^1 \frac{\sqrt{v}(1-v^n)}{1-v} dv$$

or equivalently

$$\psi\left(n+\frac{1}{2}\right) + \gamma + 2\log 2 = \frac{2}{2n+1} - \int_0^1 \frac{\sqrt{v}(1-v^n)}{1-v} dv$$

We obtain the summation

$$\sum_{n=0}^\infty \frac{\psi\left(n+\frac{1}{2}\right)}{(2n+1)^2} + (\gamma+2\log 2)\sum_{n=0}^\infty \frac{1}{(2n+1)^2} = 2\sum_{n=0}^\infty \frac{1}{(2n+1)^3} - \sum_{n=0}^\infty \frac{1}{(2n+1)^2} \int_0^1 \frac{\sqrt{v}(1-v^n)}{1-v} dv$$

Since $\varsigma(s) = \frac{1}{1-2^{-s}} \sum_{n=0}^{\infty} \frac{1}{(2n+1)^s}$ we have

$$\frac{3}{4}\varsigma(2) = \sum_{n=0}^{\infty} \frac{1}{(2n+1)^2} \qquad \frac{7}{8}\varsigma(3) = \sum_{n=0}^{\infty} \frac{1}{(2n+1)^3}$$

and hence we have

$$\sum_{n=0}^{\infty} \frac{\psi\left(n+\frac{1}{2}\right)}{(2n+1)^2} + \frac{3}{4}(\gamma + 2\log 2)\varsigma(2) = \frac{7}{4}\varsigma(3) - \int_0^1 \frac{\sqrt{v}}{1-v}\left[\frac{7}{8}\varsigma(3) - \sum_{n=0}^{\infty} \frac{v^n}{(2n+1)^2}\right] dv$$

Since

$$\sum_{n=0}^{\infty} \frac{v^n}{(2n+1)^2} = \frac{1}{4}\sum_{n=0}^{\infty} \frac{v^n}{(n+1/2)^2} = \frac{1}{4}\Phi\left(v, 2, \frac{1}{2}\right)$$

where $\Phi(v,s,a)$ is the Hurwitz-Lerch zeta function [46, p.121] defined by

$$\Phi(v,s,a) = \sum_{n=0}^{\infty} \frac{v^n}{(n+a)^s}$$

we may write this as

$$\sum_{n=0}^{\infty} \frac{\psi\left(n+\frac{1}{2}\right)}{(2n+1)^2} + \frac{3}{4}(\gamma + 2\log 2)\varsigma(2) = \frac{7}{4}\varsigma(3) - \int_0^1 \frac{\sqrt{v}}{1-v}\left[\frac{7}{8}\varsigma(3) - \frac{1}{4}\Phi\left(v, 2, \frac{1}{2}\right)\right] dv$$

In a recent preprint [10i], entitled "Determination of some generalised Euler sums involving the digamma function", using Kummer's Fourier series expansion for $\log \Gamma(x)$ it is shown that

(48.1) $$\sum_{n=0}^{\infty} \frac{\psi\left(n+\frac{1}{2}\right)}{(2n+1)^2} = -\frac{1}{8}[\gamma\pi^2 + 7\varsigma(3)]$$

and we therefore obtain the integral

(48.2) $$\int_0^1 \frac{\sqrt{v}}{1-v}\left[\frac{7}{8}\varsigma(3) - \frac{1}{4}\Phi\left(v, 2, \frac{1}{2}\right)\right] dv = \frac{21}{8}\varsigma(3) - \frac{3}{2}\log 2\varsigma(2)$$

It is also shown in [10i] that

$$(48.3) \qquad \sum_{n=0}^{\infty} \frac{\psi\left(n+\frac{1}{2}\right)}{(2n+1)^4} = -\frac{1}{96}\left[3\pi^2 \varsigma(3) + \pi^4 \gamma + 93\varsigma(5)\right]$$

Similar series involving the digamma function were extensively explored by Coffey [7i] in 2005.

**REFERENCES**

[1] V.S.Adamchik, On Stirling Numbers and Euler Sums.
J. Comput. Appl. Math.79, 119-130, 1997.
http://www-2.cs.cmu.edu/~adamchik/articles/stirling.htm

[2] V.S.Adamchik, Certain Series Associated with Catalan's Constant. Journal for Analysis and its Applications (ZAA), 21, 3 (2002), 817-826.
http://www-2.cs.cmu.edu/~adamchik/articles/csum.html

[3] P. Amore, Convergence acceleration of series through a variational approach.
arXiv:math-ph/0408036 [ps, pdf, other] 2004.

[4] G.E. Andrews, R. Askey and R. Roy, Special Functions.
Cambridge University Press, Cambridge, 1999.

[5] J. Anglesio, A fairly general family of integrals.
Amer. Math. Monthly, 104, 665-666, 1997.

[6] E.T. Bell, Exponential polynomials. Ann. of Math., 35 (1934), 258-277.

[6i] B.C. Berndt, Ramanujan's Notebooks. Parts I-III, Springer-Verlag, 1985-1991.

[7] J.M. Borwein and D.M. Bradley, Thirty-two Goldbach Variations.
Int. J. Number Theory, 2 1 (2006), 65–103. math.NT/0502034 [abs, ps, pdf, other]

[7i] M.W. Coffey, On one-dimensional digamma and polygamma series related to the evaluation of Feynman diagrams.
J. Comput. Appl. Math, 183, 84-100, 2005. math-ph/0505051 [abs, ps, pdf, other]

[8] M.W. Coffey, A set of identities for a class of alternating binomial sums arising in computing applications. 2006.
arXiv:math-ph/0608049v1

[9] C.B. Collins, The role of Bell polynomials in integration.
J. Comput. Appl. Math. 131 (2001) 195-211.

[10] L. Comtet, Advanced Combinatorics, Reidel, Dordrecht, 1974.

[10i] D.F. Connon, Determination of some generalised Euler sums involving the digamma function. 2008. arXiv:0802.1440 [pdf]

[11] D.F. Connon, Some series and integrals involving the Riemann zeta function, binomial coefficients and the harmonic numbers. Volume I, 2007. arXiv:0710.4022 [pdf]

[12] D.F. Connon, Some series and integrals involving the Riemann zeta function, binomial coefficients and the harmonic numbers. Volume II(a), 2007. arXiv:0710.4023 [pdf]

[13] D.F. Connon, Some series and integrals involving the Riemann zeta function, binomial coefficients and the harmonic numbers. Volume II(b), 2007. arXiv:0710.4024 [pdf]

[14] D.F. Connon, Some series and integrals involving the Riemann zeta function, binomial coefficients and the harmonic numbers. Volume III, 2007. arXiv:0710.4025 [pdf]

[15] D.F. Connon, Some series and integrals involving the Riemann zeta function, binomial coefficients and the harmonic numbers. Volume IV, 2007. arXiv:0710.4028 [pdf]

[16] D.F. Connon, Some series and integrals involving the Riemann zeta function, binomial coefficients and the harmonic numbers. Volume V, 2007. arXiv:0710.4047 [pdf]

[17] D.F. Connon, Some series and integrals involving the Riemann zeta function, binomial coefficients and the harmonic numbers. Volume VI, 2007. arXiv:0710.4032 [pdf]

[18] M.A. Coppo, La formule d'Hermite revisitée. 2003. http://math1.unice.fr/~coppo/

[19] P. Flajolet and R. Sedgewick, Mellin Transforms and Asymptotics: Finite Differences and Rice's Integrals. Theor. Comput. Sci. 144, 101-124, 1995. Mellin Transforms and Asymptotics : Finite Differences and Rice's Integrals (117kb),

[20] J. Fleisher, A.V. Kotikov and O.L. Veretin, Analytic two-loop results for self energy- and vertex-type diagrams with one non-zero mass. hep-ph/9808242 [abs, ps, pdf, other] Nucl. Phys. B547 (1999) 343-374.

[21] H.W. Gould, Some Relations involving the Finite Harmonic Series. Math. Mag., 34, 317-321, 1961.

[22] H.W. Gould, Combinatorial Identities.Rev.Ed.University of West
     Virginia, U.S.A., 1972.

[23] H.W. Gould, Explicit formulas of Bernoulli Numbers.
      Amer. Math. Monthly, 79, 44-51, 1972.

[24] R.L. Graham, D.E. Knuth and O. Patashnik, Concrete Mathematics. Second Ed.
     Addison-Wesley Publishing Company, Reading, Massachusetts, 1994.

[25] H. Hasse, Ein Summierungsverfahren für Die Riemannsche $\varsigma$ - Reithe.
     Math. Z.32, 458-464, 1930.
     http://dz-srv1.sub.uni-goettingen.de/sub/digbib/loader?ht=VIEW&did=D23956&p=462

[26] W.P. Johnson, The Curious History of Faà di Bruno's Formula.
     Amer. Math. Monthly 109,217-234, 2002.

[27] P. Kirschenhofer, A Note on Alternating Sums. The Electronic Journal of
     Combinatorics 3 (2), #R7, 1996. R7: Peter Kirschenhofer

[28] K. Knopp, Theory and Application of Infinite Series. Second English Edition.
     Dover Publications Inc, New York, 1990.

[29] K.S. Kölbig and W. Strampp Some infinite integrals with powers of logarithms and
     the complete Bell polynomials.
     J. Comput. Appl. Math. 69 (1996) 39-47.
     Also available electronically at:
     An integral by recurrence and the Bell polynomials.
     CERN/Computing and Networks Division, CN/93/7, 1993.
     http://cdsweb.cern.ch/record/249027/

[30] K.S. Kölbig ,The complete Bell polynomials for certain arguments in terms of
     Stirling numbers of the first kind.
     J. Comput. Appl. Math. 51 (1994) 113-116.
     Also available electronically at:
     A relation between the Bell polynomials at certain arguments and a Pochhammer
     symbol. CERN/Computing and Networks Division, CN/93/2, 1993.
     http://doc.cern.ch/archive/electronic/other/preprints//CM-P/CM-P00065731.pdf

[30i] K.S. Kölbig, The polygamma function $\psi^{(k)}(x)$ for $x = \frac{1}{4}$ and $x = \frac{3}{4}$,
      J. Comput. Appl. Math. 75 (1996) 43-46.

[31] P.J. Larcombe, E.J. Fennessey and W.A. Koepf, Integral proofs of two
     alternating sign binomial coefficient identities.
     Util. Math.66, 93-103 (2004) http://citeseer.ist.psu.edu/598454.html

[32] M.E. Levenson, J.F. Locke and H. Tate, Amer. Math. Monthly, 45, 56-58, 1938.

[33] Z.R. Melzak, Companion to Concrete Mathematics. Wiley-Interscience, New York, 1973.

[34] N. Nielsen, Die Gammafunktion. Chelsea Publishing Company, Bronx and New York, 1965.

[35] N.E. Nörlund, Vorlesungen über Differenzenrechnung. Chelsea, 1954.
http://dz-srv1.sub.uni-goettingen.de/cache/browse/AuthorMathematicaMonograph,WorkContainedN1.html

[36] N.E. Nörlund, Leçons sur les séries d'interpolation.
Paris, Gauthier-Villars, 1926.

[37] G. Póyla and G. Szegö, Problems and Theorems in Analysis, Vol.I
Springer-Verlag, New York 1972.

[38] R.K. Raina and R.K. Ladda, A new family of functional series relations involving digamma functions.
Ann. Math. Blaise Pascal, Vol. 3, No. 2, 1996, 189-198.
http://www.numdam.org/item?id=AMBP_1996__3_2_189_0

[39] Srinivasa Ramanujan, Notebooks of Srinivasa Ramanujan, Vol.1, Tata Institute of Fundamental Research, Bombay, 1957.

[40] H. Ruben, A Note on the Trigamma Function.
Amer. Math. Monthly, 83, 622-623, 1976.

[41] L.-C. Shen, Remarks on some integrals and series involving the Stirling numbers and $\varsigma(n)$. Trans. Amer. Math. Soc. 347, 1391-1399, 1995.

[42] A. Snowden, Collection of Mathematical Articles. 2003.
http://www.math.princeton.edu/~asnowden/math-cont/dorfman.pdf

[43] J. Sondow, Analytic Continuation of Riemann's Zeta Function and Values at Negative Integers via Euler's Transformation of Series. Proc.Amer.Math.Soc. 120,421-424, 1994.
http://home.earthlink.net/~jsondow/id5.html

[44] J. Spieß, Some identities involving harmonic numbers.
Math. of Computation, 55, No.192, 839-863, 1990.

[45] W.G. Spohn; A.S. Adikesavan; H.W.Gould.
Amer. Math. Monthly, 75, 204- 205,1968.

[46] H.M. Srivastava and J. Choi, Series Associated with the Zeta and Related Functions. Kluwer Academic Publishers, Dordrecht, the Netherlands, 2001.

[47] The Mactutor History of Mathematics archive.

http://www-history.mcs.st-andrews.ac.uk/Mathematicians/Faa_di_Bruno.html

[48] J.A.M. Vermaseren, Harmonic sums, Mellin transforms and integrals.
Int. J. Mod. Phys. A14 (1999) 2037-2076
http://arXiv.org/abs/hep-ph/9806280

[49] E.T. Whittaker and G.N. Watson, A Course of Modern Analysis: An Introduction to the General Theory of Infinite Processes and of Analytic Functions; With an Account of the Principal Transcendental Functions. Fourth Ed., Cambridge University Press, Cambridge, London and New York, 1963.

[50] H.S. Wilf, The asymptotic behaviour of the Stirling numbers of the first kind.
Journal of Combinatorial Theory Series A, 64, 344-349, 1993.
http://citeseer.ist.psu.edu/520553.html

[51] De-Yin Zheng, Further summation formulae related to generalized harmonic numbers. J. Math. Anal. Appl. 335 (2007) 692-706.

Donal F. Connon
Elmhurst
Dundle Road
Matfield, Kent TN12 7HD
dconnon@btopenworld.com